\documentclass[11pt,a4paper]{article}
\usepackage{amsmath, amssymb}
\usepackage{latexsym, color}
\usepackage{epsfig}

\numberwithin{equation}{section}

\begin{document}

\newtheorem{thm}{Theorem}[section]
\newtheorem{cor}[thm]{Corollary}
\newtheorem{lem}[thm]{Lemma}
\newtheorem{prop}[thm]{Proposition}
\newtheorem{definition}[thm]{Definition}
\newtheorem{rem}[thm]{Remark}
\newtheorem{Ex}[thm]{EXAMPLE}
\def\nm{\noalign{\medskip}}

\bibliographystyle{plain}

%\numberwithin{equation}{section}

\newcommand{\qed}{\hfill \ensuremath{\square}}
\newcommand{\ds}{\displaystyle}
\newcommand{\pf}{\noindent {\sl Proof}. \ }
\newcommand{\p}{\partial}
\newcommand{\pd}[2]{\frac {\p #1}{\p #2}}
\newcommand{\norm}[1]{\| #1 \|}
\newcommand{\dbar}{\overline \p}
\newcommand{\eqnref}[1]{(\ref {#1})}
\newcommand{\na}{\nabla}
\newcommand{\one}[1]{#1^{(1)}}
\newcommand{\two}[1]{#1^{(2)}}

\newcommand{\Abb}{\mathbb{A}}
\newcommand{\Cbb}{\mathbb{C}}
\newcommand{\Ibb}{\mathbb{I}}
\newcommand{\Nbb}{\mathbb{N}}
\newcommand{\Kbb}{\mathbb{K}}
\newcommand{\Rbb}{\mathbb{R}}
\newcommand{\Sbb}{\mathbb{S}}

\renewcommand{\div}{\mbox{div}~}

\newcommand{\la}{\langle}
\newcommand{\ra}{\rangle}

\newcommand{\Hcal}{\mathcal{H}}
\newcommand{\Lcal}{\mathcal{L}}
\newcommand{\Kcal}{\mathcal{K}}
\newcommand{\Dcal}{\mathcal{D}}
\newcommand{\Pcal}{\mathcal{P}}
\newcommand{\Qcal}{\mathcal{Q}}
\newcommand{\Scal}{\mathcal{S}}

%%%%%%%%%%%%%%%%%%%%%%%%%%%%%%%%%%%%%%%%%%%%
% define bold face
%%%%%%%%%%%%%%%%%%%%%%%%%%%%%%%%%%%%%%%%%%%%
\def\Ba{{\bf a}}
\def\Bb{{\bf b}}
\def\Bc{{\bf c}}
\def\Bd{{\bf d}}
\def\Be{{\bf e}}
\def\Bf{{\bf f}}
\def\Bg{{\bf g}}
\def\Bh{{\bf h}}
\def\Bi{{\bf i}}
\def\Bj{{\bf j}}
\def\Bk{{\bf k}}
\def\Bl{{\bf l}}
\def\Bm{{\bf m}}
\def\Bn{{\bf n}}
\def\Bo{{\bf o}}
\def\Bp{{\bf p}}
\def\Bq{{\bf q}}
\def\Br{{\bf r}}
\def\Bs{{\bf s}}
\def\Bt{{\bf t}}
\def\Bu{{\bf u}}
\def\Bv{{\bf v}}
\def\Bw{{\bf w}}
\def\Bx{{\bf x}}
\def\By{{\bf y}}
\def\Bz{{\bf z}}
\def\BA{{\bf A}}
\def\BB{{\bf B}}
\def\BC{{\bf C}}
\def\BD{{\bf D}}
\def\BE{{\bf E}}
\def\BF{{\bf F}}
\def\BG{{\bf G}}
\def\BH{{\bf H}}
\def\BI{{\bf I}}
\def\BJ{{\bf J}}
\def\BK{{\bf K}}
\def\BL{{\bf L}}
\def\BM{{\bf M}}
\def\BN{{\bf N}}
\def\BO{{\bf O}}
\def\BP{{\bf P}}
\def\BQ{{\bf Q}}
\def\BR{{\bf R}}
\def\BS{{\bf S}}
\def\BT{{\bf T}}
\def\BU{{\bf U}}
\def\BV{{\bf V}}
\def\BW{{\bf W}}
\def\BX{{\bf X}}
\def\BY{{\bf Y}}
\def\BZ{{\bf Z}}

%%%%%%%%%%%%%%%%%%%%%%%%%%%%%%%%%%%%%%%%%%%%%%%%%%%%
% Abbreviate definitions of greek symbols
%%%%%%%%%%%%%%%%%%%%%%%%%%%%%%%%%%%%%%%%%%%%%%%%%%%%

\newcommand{\Ga}{\alpha}
\newcommand{\Gb}{\beta}
\newcommand{\Gd}{\delta}
\newcommand{\Ge}{\epsilon}
\newcommand{\Gve}{\varepsilon}
\newcommand{\Gf}{\phi}
\newcommand{\Gvf}{\varphi}
\newcommand{\Gg}{\gamma}
\newcommand{\Gc}{\chi}
\newcommand{\Gi}{\iota}
\newcommand{\Gk}{\kappa}
\newcommand{\Gvk}{\varkappa}
\newcommand{\Gl}{\lambda}
\newcommand{\Gn}{\eta}
\newcommand{\Gm}{\mu}
\newcommand{\Gv}{\nu}
\newcommand{\Gp}{\pi}
\newcommand{\Gt}{\theta}
\newcommand{\Gvt}{\vartheta}
\newcommand{\Gr}{\rho}
\newcommand{\Gvr}{\varrho}
\newcommand{\Gs}{\sigma}
\newcommand{\Gvs}{\varsigma}
\newcommand{\Gj}{\tau}
\newcommand{\Gu}{\upsilon}
\newcommand{\Go}{\omega}
\newcommand{\Gx}{\xi}
\newcommand{\Gy}{\psi}
\newcommand{\Gz}{\zeta}
\newcommand{\GD}{\Delta}
\newcommand{\GF}{\Phi}
\newcommand{\GG}{\Gamma}
\newcommand{\GL}{\Lambda}
\newcommand{\GP}{\Pi}
\newcommand{\GT}{\Theta}
\newcommand{\GS}{\Sigma}
\newcommand{\GU}{\Upsilon}
\newcommand{\GO}{\Omega}
\newcommand{\GX}{\Xi}
\newcommand{\GY}{\Psi}

%%%%%%%%%%
\newcommand{\beq}{\begin{equation}}
\newcommand{\eeq}{\end{equation}}

\title{Spectral analysis of the
Neumann-Poincar\'e operator and characterization of the gradient blow-up\thanks{\footnotesize This work was
supported by the ERC Advanced Grant Project MULTIMOD--267184 and
NRF grants No. 2009-0085987, 2010-0004091, and 2010-0017532, and by Hankuk University
of Foreign Studies Research Fund of 2012}}

\author{Habib Ammari\thanks{\footnotesize Department of Mathematics and Applications, Ecole Normale Sup\'erieure,
45 Rue d'Ulm, 75005 Paris, France (habib.ammari@ens.fr).} \and
Giulio Ciraolo\thanks{\footnotesize  Dipartimento di Matematica e
Informatica, Universit\`a di Palermo, Via Archirafi 34, 90123,
Palermo,  Italy
  (g.ciraolo@math.unipa.it).} \and Hyeonbae
Kang\thanks{Department of Mathematics, Inha University, Incheon
402-751, Korea (hbkang@inha.ac.kr, hdlee@inha.ac.kr).}  \and
Hyundae Lee\footnotemark[4]  \and KiHyun Yun\thanks{\footnotesize
Department of Mathematics, Hankuk University of Foreign Studies,
Youngin-si, Gyeonggi-do 449-791, Korea
(gundam@hufs.ac.kr).}}
%\date{}

\maketitle

\begin{abstract}
When perfectly conducting or insulating inclusions are closely
located, stress which is the gradient of the solution to the
conductivity equation can be arbitrarily large as the distance
between two inclusions tends to zero. It is important to precisely
characterize the blow-up of the gradient. In this paper we show
that the blow-up of the gradient can be characterized by a
singular function defined by the single layer potential of an
eigenfunction corresponding to the eigenvalue ${1}/{2}$ of a
Neumann-Poincar\'e type operator defined on the boundaries of the
inclusions. By comparing the singular function with the one
corresponding to two disks osculating to the inclusions, we
quantitatively characterize the blow-up of the gradient in terms
of explicit functions.
\end{abstract}

\noindent {\footnotesize {\bf Mathematics subject classification
(MSC2000)}: 35J25, 73C40}


%%%%%%%%%%%%%%%%%%%%%%%%%%%%%%%%%%%
\section{Introduction}
%%%%%%%%%%%%%%%%%%%%%%%%%%%%%%%%%%%

Let $D_1$ and $D_2$ be bounded simply connected
domains in $\Rbb^d$, $d=2,3$, whose boundary regularity will be specified later. Suppose that they are conductors, whose conductivity is $k$, $0 < k \neq 1 < \infty$, embedded in the background with conductivity 1. Let $\Gs$ denote the
conductivity distribution, {\it i.e.},
\beq\label{conddis}
\Gs= k \chi(D_1 \cup D_2) + \chi(\Rbb^d \setminus (D_1 \cup D_2)),
\eeq
where $\chi$ is the characteristic function. We consider the following elliptic problem: for a given entire
harmonic function $h$ in $\Rbb^d$,
 \beq\label{eqnk}
 \begin{cases}
 \ds \nabla  \cdot \Gs \nabla u =0 \quad \mbox{in } \Rbb^d, \\
 u(\Bx)- h(\Bx) = O(|\Bx|^{1-d}) \quad \mbox{as } |\Bx| \to
 \infty.
 \end{cases}
 \eeq
Let
\beq
\Ge:= \mbox{dist}(D_1, D_2),
\eeq
and assume that $\Ge$ is small. We emphasize that the shapes of $D_1$ and $D_2$ do not depend on $\Ge$. More precisely, there are fixed domains $\widetilde{D}_1$ and $\widetilde{D}_2$ such that $D_j$ is a translate of $\widetilde{D}_j$, namely, there are vectors $\Ba_1$ and $\Ba_2$ such that
 \beq\label{ajtrans}
 D_j=\widetilde{D}_j + \Ba_j, \quad j=1,2.
 \eeq
The problem is to estimate $|\nabla u|$ in terms of $\Ge$ when $\Ge$ tends to $0$, or to characterize the asymptotic singular behavior of $\nabla u$ as $\Ge \to 0$.

If $k$ stays away from $0$ and $\infty$, {\it i.e.}, $c_1 < k < c_2$ for some positive constants
$c_1$ and $c_2$, then $|\nabla u|$ is bounded regardless of $\Ge$ as was proved
in \cite{BV, LV, LN}. In fact, it is proved that the $\mathcal{C}^{1,\Ga}$ norm of $u$ is bounded regardless of $\Ge$ when $\p D_1$ and $\p D_2$ are $\mathcal{C}^{2,\Ga}$ smooth. However, if $k$ degenerates to either $\infty$ (perfectly conducting case) or $0$ (insulating case), the ellipticity holds only outside $D_1$ and $D_2$ and completely different phenomena occur.

When $k=\infty$, the problem becomes
\beq\label{inftycond}
\quad \left\{
\begin{array}{ll}
\ds \Delta u  = 0 \quad\mbox{in } \Rbb^d \setminus \overline{(D_1 \cup D_2)},\\
\ds u = \Gl_i \ (\mbox{constant}) \quad\mbox{on } \p D_i,~i=1,2, \\
\ds u(\Bx) - h(\Bx) =O(|\Bx|^{1-d}) \quad \mbox{as }
|\Bx|\rightarrow \infty,
\end{array}
\right. \eeq where the constants $\Gl_i$ are determined by the
conditions
 \beq
 \int_{\p D_1} \frac {\p u}{\p \nu^{(1)}} \Big|_+  = \int_{\p D_2} \frac {\p u}{\p \nu^{(2)}} \Big|_+ =0 ,
 \eeq
with $\nu^{(j)}$ being the outward unit normal to $\partial D_j$,
$j=1,2$. In that case, $\nabla u$ may blow up as $\Ge$ tends to $0$.

As shown in  \cite{keller, BC, AKL, AKLLL, Y, Y2, BLY}, in two
dimensions the generic rate of gradient blow-up is $\Ge^{-1/2}$,
while it is $|\Ge \log \Ge|^{-1}$ in three dimensions \cite{BLY,
BLY2, LY}. The blow-up of the gradient may or may not occur
depending on the background potential (the harmonic function $h$
in \eqnref{eqnk}) and those background potentials which actually
make the gradient blow up are characterized in \cite{AKLLZ} when $D_1$ and $D_2$ are disks. In
two dimensions, the perfectly insulating case, where $k=0$, can be
dealt with using the conjugate relation (see \cite{keller, AKL})
and in this case the blow-up rate is also $\Ge^{-1/2}$. It
is a challenging open problem to clarify whether $|\nabla u|$ may
blow up or not in the insulating case in three dimensions and to
find the blow-up rate if the blow-up occurs. It is also a quite
interesting problem to clarify the dependence of $|\nabla u|$ on $k$
as $k \to \infty$ or $k \to 0$. In this relation we mention that a
precise dependence on $k$ when $D_1$ and $D_2$ are disks was
shown in \cite{AKL, AKLLL}. It is worth mentioning that a
similar blow-up phenomenon for the $p$-Laplacian equation was
investigated in \cite{GN}.

Above mentioned results for $k=\infty$ are estimates of $|\nabla
u|$ in upper and lower bounds sense, namely,
\beq
C_1 \psi(\Ge) \le \| \nabla u \|_\infty \le C_2 \psi(\Ge)
\eeq
for some implicit constants
$C_1$ and $C_2$ where $\psi(\Ge)= \Ge^{-1/2}$ in two dimensions
and $\psi(\Ge)= |\Ge \log \Ge|^{-1}$ in three dimensions. The
constant $C_1$ can be zero or positive depending on the background
potential $h$. In order to have a better understanding of the
concentration of the gradient it is desirable to pursue deeper
investigation on the blow-up nature of $\nabla u$. In this
direction there is a recent work \cite{KLY} where the blow-up
nature of $\nabla u$ is characterized in terms of an explicit
singular function. It is shown that if $D_1=B_1$ and $D_2=B_2$ are
circular inclusions of radius $r_1$ and $r_2$, respectively, and
$k=\infty$, then
 \beq\label{uahb}
 u(\Bx)= \frac{2 r_1 r_2}{r_1 + r_2} (\Bn \cdot \nabla h) (\Bp)
 \left( \ln |\Bx-\Bp_1|- \ln |\Bx-\Bp_2|\right) + r(\Bx),
 \eeq
for $\Bx \in \Rbb^2 \setminus (B_1 \cup B_2)$, where $\Bp_1 \in D_1$ is the fixed point of $R_1 R_2$ where $R_j$
is the reflection with respect to $\p B_j$, $j=1,2$, $\Bp_2
\in B_2$ is the fixed point of $R_2 R_1$, $\Bn$ is the
 unit vector in the direction of $\Bp_2 -\Bp_1$, and $\Bp$ is the middle
point of the shortest line segment connecting $\p B_1$ and $\p
B_2$. In \eqnref{uahb},
$\nabla r$ is bounded independently of $\Ge$ and thus the blow-up
of $\nabla u$ is completely characterized by the singular function
 \beq\label{singular}
 q_B(\Bx):= \frac{1}{2\pi} \left( \ln |\Bx-\Bp_1|- \ln |\Bx-\Bp_2|\right).
 \eeq

The purpose of this paper is to establish a characterization of
the blow-up of $\nabla u$ similar to \eqnref{uahb} when $D_1$ and
$D_2$ are strictly convex simply connected domains in $\Rbb^2$. In
doing so, the Neumann-Poincar\'e (NP) operator denoted by $\Kbb^*$
and defined on $L^2(\p D_1) \times L^2(\p D_2)$ plays a crucial
role. The NP operator is a classical notion and  appears naturally
when we solve the boundary value problems using layer potentials.
It also appears naturally when we solve the transmission problem
\eqnref{eqnk}. See the next section for a definition of the NP
operator for the problem of this paper. This operator has $1/2$ as
an eigenvalue with multiplicity $2$. (If the inclusion has $N$
simply connected components, then the multiplicity of $1/2$ is
$N$.) If two inclusions are disks, then $(\pd{q_B}{\nu^{(1)}}|_{\p
D_1}, \pd{q_B}{\nu^{(2)}}|_{\p D_2})^T$, where $q_B$ is given by
\eqnref{singular}, is an eigenfunction of the NP operator on
$L^2(\p D_1) \times L^2(\p D_2)$ corresponding to $1/2$. (This
fact was also observed in \cite{BT}.) Here, $T$ denotes the
transpose.

Let $\Bg=(\one{g},\two{g})^T$ be the eigenfunction of $\Kbb^*$ on
$L^2(\p D_1) \times L^2(\p D_2)$ corresponding to the eigenvalue
$1/2$ and satisfying \beq\label{plusone} \int_{\p D_1} \one{g}
d\Gs = 1, \quad \int_{\p D_2} \two{g} d\Gs =-1. \eeq We will prove
that such an eigenfunction does exist. Let $h$ be the background
harmonic function introduced in \eqnref{eqnk} and let \beq \Bh:=
\begin{bmatrix} h|_{\p D_1} \\ h|_{\p D_2} \end{bmatrix}. \eeq Let
$\la \Bh, \Bg \ra$ be the inner product on $L^2(\p D_1) \times
L^2(\p D_2)$, {\it i.e.},
\beq
\la \Bh, \Bg \ra = \int_{\p D_1} h
\one{g} d\Gs + \int_{\p D_2} h \two{g} d\Gs .
\eeq
With these notions in hand, we can state the main result of this paper.
\begin{thm}\label{thm3}
Let $D_1$ and $D_2$ be strictly convex simply connected domains in $\Rbb^2$ with $\mathcal{C}^{2,\alpha}$ smooth
boundaries for some $\alpha \in (0,1]$. Let $\Bz_1 \in \p D_1$ and $\Bz_2 \in \p D_2$ be the closest points, and let
$\Ge:= \mbox{dist}(D_1, D_2)=|\Bz_1-\Bz_2|$, $\Gk_j$ be the curvature of $\p D_j$ at $\Bz_j$, $B_j$
be the disk osculating to $D_j$ at $\Bz_j$, $j=1,2$, and $q_B$ be the singular function in \eqnref{singular}
associated with disks $B_1$ and $B_2$.  Then, the solution $u$ to \eqnref{inftycond} satisfies
\beq\label{charcaterize}
u(\Bx) = -\frac{\sqrt 2 \pi \la \Bh, \Bg \ra}{\sqrt{\Ge(\Gk_1+\Gk_2)}} \Ga_\Ge q_B (\Bx) + r (\Bx),
\quad \Bx \in \Rbb^2 \setminus (D_1 \cup D_2),
\eeq
where $\Ga_\Ge$ is a constant bounded independently of $\Ge$ of the form
\beq
\quad\Ga_\Ge =  \left\{
\begin{array}{ll}
\ds
1+ O ({\Ge}^{ {\alpha}/ 2 } ) \quad &\mbox{if  } \, \alpha \in (0,1),\\
1+ O (|\sqrt{\Ge} \ln\Ge|) \quad &\mbox{if  } \, \alpha =1,
\end{array}
\right.
\quad\mbox{as } \Ge \to 0,
\eeq
and
\beq
\norm{\nabla r}_{L^{\infty} (\GO)} \leq C
\eeq
for some constant
$C$ independent of $\Ge$. Here $\GO = \GO_1 \setminus (D_1 \cup
D_2)$ and $\GO_1$ is an open set containing $\overline{D_1 \cup
D_2}$.
\end{thm}

We emphasize that \eqnref{charcaterize} is a pointwise relation
and hence describes the behavior of $\nabla u(\Bx)$ in terms of
the gradient of the function $q_B(\Bx)$. One can see from the explicit expression \eqnref{singular} that $|\nabla q_B|$ attains its maximum at $\Bz_1$ and $\Bz_2$, and that
 \beq\label{qBmax2}
 |\nabla q_B (\Bz_j)| = \frac{\sqrt{\Gk_1+\Gk_2}}{\sqrt 2 \pi} \frac{1}{\sqrt{\Ge}} + O(1).
 \eeq
(See \eqnref{pjest}.) So, \eqnref{charcaterize} shows that $|\nabla u|$ is bounded regardless of $\Ge$ if $\la \Bh, \Bg \ra = 0$. Moreover, it yields a new improved estimate:
 \beq\label{linfty}
 \| \nabla u \|_{L^\infty(\Rbb^2 \setminus (D_1 \cup D_2))} = \frac{\Ga_\Ge |\la \Bh, \Bg \ra|}{\Ge} + O(1), \quad\mbox{as } \Ge \to 0.
 \eeq
Since
\beq
|\la \Bh, \Bg \ra| \le C \sqrt{\Ge}
\eeq
for some constant $C$ independent of $\Ge$ as shown in \eqnref{hgbound}, we can also infer from \eqnref{linfty} that the generic rate of blow-up is $\Ge^{-1/2}$.

The (global) strict convexity assumption of $D_1$ and $D_2$ in Theorem \ref{thm3} can be relaxed a little. Instead, if we let $\Bz=\frac{\Bz_1+\Bz_2}{2}$, then it is enough to assume that there are $\Gd$ and $c_0$ such that $\p D_j \cap B_\Gd(\Bz)$ is strictly convex for $j=1,2$, and
\beq
\mbox{dist} (\p D_1 \setminus B_\Gd(\Bz), \p D_2 \setminus B_\Gd(\Bz)) \ge c_0.
\eeq
This can be shown by exactly the same proofs as in this paper.

We also obtain similar results for the insulating case and
boundary value problems. The problem for the insulating case,
obtained by taking the limit as $k\rightarrow 0$ of \eqnref{eqnk},
is given by
\beq\label{zeocond}
\quad \left\{
\begin{array}{ll}
\ds \Delta u  = 0 \quad\mbox{in } \Rbb^d \setminus \overline{(D_1 \cup D_2)},\\
\ds \pd{u}{\nu^{(i)}} \Big|_{+}=0 \quad\mbox{on } \p D_i,~i=1,2, \\
\ds u(\Bx) - h(\Bx) =O(|\Bx|^{1-d}) \quad \mbox{as } |\Bx|\rightarrow \infty.
\end{array}
\right. \eeq  If $u$ is the solution to \eqnref{zeocond}, then its
harmonic conjugate $u^\perp$ is the solution to \eqnref{inftycond}
with $h$ replaced with its harmonic conjugate $h^\perp$. Therefore, we may apply Theorem \ref{thm3} to $u^\perp$ to obtain an analogous result for $u$.

On the other hand, if we consider the boundary value problem
 \beq
 \nabla  \cdot \Gs \nabla u =0 \quad \mbox{in } \GO ,
 \eeq
with either Dirichlet or Neumann boundary conditions on the
boundary $\p\GO$ of a smooth domain $\GO$ containing $D_1$ and
$D_2$ and $k \rightarrow 0$ or $\infty$, then, in this case, the
harmonic function $h$ defined by \beq\label{hbvp} h(\Bx) = -
\frac{1}{2\pi} \int_{\p \GO} \ln |\Bx-\By| \pd{u}{\nu}(\By) \,
d\Gs(\By) -
 \frac{1}{2\pi} \int_{\p \GO} \frac{\langle \Bx-\By, \nu(\By) \rangle}{|\Bx-\By|^2} u(\By) \,
 d\Gs(\By),
\eeq with $\nu$ being the outward normal to $\partial \Omega$,
plays the role of $h$ in the whole space problem
\eqnref{inftycond} or \eqnref{zeocond}. A similar result on the
characterization of the gradient blow-up can be obtained using
exactly the same arguments as for the whole space problem.

The results of this paper can be applied for solving two
longstanding problems. The first one is the study of material
failure. In fact, the problem of estimation of the gradient
blow-up was raised by Babu\u{s}ka in relation to the study of
material failure of composites \cite{bab}. In composites which
consist of inclusions and the matrix, some inclusions may be
closely located and stress occurs in between them. The problems
\eqnref{eqnk}, \eqnref{inftycond} and \eqnref{zeocond} are
conductivity or anti-plane elasticity equation, and
$\nabla u$ represents the shear stress tensor. So results of this paper provide clear quantitative
understanding of the stress concentration, which will be a
fundamental ingredient in the study of material failure.

The second application is computation of the electrical field in the presence of closely located inclusions
with extreme conductivities ($0$ or $\infty$) which is known to be a hard problem. Because $|\nabla u|$ becomes arbitrarily large, we need fine
meshes to compute $\nabla u$ numerically. Since \eqnref{charcaterize} for example provides complete characterization
of the singular behavior of $\nabla u$, the complexity of computation can be greatly reduced by
removing the singular term there. In fact, effectiveness of this scheme is already demonstrated
in \cite{KLY} when inclusions are disks using \eqnref{uahb}. It is worth mentioning that unlike \eqnref{uahb} where the coefficient of $q_B$ is explicitly determined by $h$, computation of the
constant $\Ga_\Ge \la \Bh, \Bg \ra$ in \eqnref{charcaterize} may cause a problem when inclusions are of general shape.
 We will investigate this issue in a forthcoming work.

This paper is organized as follows. In the next section we
introduce the single layer potential and define the
Neumann-Poincar\'e operator. In section 3, we construct
eigenfunctions of the NP operator corresponding to the eigenvalue
$1/2$ and prove that its multiplicity is 2. In section 4, we
construct a singular function using eigenfunctions constructed in
the previous section and characterize the gradient blow-up in
terms of the singular function. In section 5, we estimate the
potential difference of the solution to \eqnref{inftycond}.
Section 6 is to prove Theorem \ref{thm3}. Sections 7 and 8 are for
the insulating case and the boundary value problem, respectively.
In the last section we prove a lemma used in Section 6.

%%%%%%%%%%%%%%%%%%%%%%%%%%%%%%%%%%%%%%%%%%%
\section{Preliminaries}
%%%%%%%%%%%%%%%%%%%%%%%%%%%%%%%%%%%%%%%%%%%

Let $D$ be a bounded simply connected domain in $\Rbb^d$, $d=2,3$,
with a Lipschitz boundary. The single layer potential
$\Scal_D[\Gvf]$ of a function $\Gvf \in L^2 (\p D)$ is defined as
 $$
 \mathcal{S}_D [\Gvf] (\Bx) = \int_{\p D} \GG(\Bx-\By) \Gvf(\By)~d\sigma (\By), \quad \Bx \in \Rbb^d,
 $$
where $\GG$ is the fundamental solution to the Laplacian, {\it i.e.},
 \beq\label{gammacond}
 \GG (\Bx) =
 \begin{cases}
 \ds \frac{1}{2\pi} \ln |\Bx|\;, \quad & d=2 \;, \\ \nm \ds
 -\frac{1}{4\pi} |\Bx|^{-1}\;, \quad & d = 3 \;.
 \end{cases}
 \eeq
Then, for $\Gvf \in L^2 (\p D)$, we have
 \beq\label{jump}
 \pd{}{\nu} \mathcal{S}_D [\Gvf] \Big|_{\pm} (\Bx) =
 \left(\pm {\frac 1 2} + \mathcal{K}_D ^* \right) [\Gvf]  (\Bx)~~\mbox{a.e. }  \Bx\in \p D,
 \eeq
where
 \beq
 \Kcal_D^* [\Gvf] (\Bx) = \int_{\p D} \frac{\p \GG(\Bx-\By)}{\p\nu(\Bx)} \Gvf(\By)~d\sigma (\By), \quad \Bx \in \p D.
 \eeq
Here, $\pd{}{\nu}$ denotes the normal derivative and the subscripts $+$ and $-$ represent the
limits from outside and inside $D$, respectively. The operator $\Kcal_D^*$ is called the Neumann-Poincar\'e
(NP) operator.

We now consider the configuration where there are two inclusions $D_1$ and $D_2$ which are closely located.
Suppose that the conductivity of the inclusions is $k \neq 1$ while that of the background is $1$, so
that the conductivity distribution is given by \eqnref{conddis}. For a given entire harmonic
function $h$ in $\Rbb^d$, we consider the problem \eqnref{eqnk}.

It is known (see for example \cite{KS96, KS2000}) that the solution $u$ to \eqnref{inftycond} can be represented as
 \beq\label{urepres}
 u (\Bx)= h (\Bx) + \Scal_{D_1}[\one{\Gvf}] (\Bx) +  \Scal_{D_2}[\two{\Gvf}] (\Bx), \quad \Bx \in \Rbb^d \setminus (D_1 \cup D_2)
 \eeq
for a pair of functions $(\one{\Gvf}, \two{\Gvf}) \in L^2_0(\p
D_1) \times L^2_0(\p D_2)$ ($L^2_0$ denotes the set of $L^2$
functions with mean zero).  Since $u$ is constant on $\p D_1$ and
$\p D_2$, we have
 $$
 \pd{}{\nu^{(j)}} \big( \Scal_{D_1}[\one{\Gvf}] +  \Scal_{D_2}[\two{\Gvf}] \big) \Big|_{-} =
 - \pd{h}{\nu^{(j)}} \quad\mbox{on } \p D_j, \ j=1,2,
 $$
which, according to \eqnref{jump}, may be written as
 \begin{align}
 \left( \frac{1}{2} I - \Kcal_{D_1}^* \right) [\one{\Gvf}] -
 \pd{}{\one{\nu}} \mathcal{S}_{D_2} [\two{\Gvf}] = \pd{h}{\one{\nu}} \quad\mbox{on } \p D_1, \\
 - \pd{}{\two{\nu}} \mathcal{S}_{D_2} [\one{\Gvf}]
 + \left( \frac{1}{2} I - \Kcal_{D_2}^* \right) [\two{\Gvf}] = \pd{h}{\two{\nu}} \quad\mbox{on } \p D_2.
 \end{align}
Here, $\pd{h}{\nu^{(j)}}$ denotes the outward normal derivative on $\p
D_j$, $j=1,2$. This system of integral equations can be written in
a condensed form as
 \beq
 \big( \frac{1}{2} \Ibb - \Kbb^* \big) [\Gvf] = \p h,
 \eeq
where
 \beq\label{Ibb}
 \Ibb = \begin{bmatrix} I  & 0 \\  0  &  I \end{bmatrix}, \quad \Kbb^* :=\begin{bmatrix}
  \Kcal_{D _1}^* & \frac {\p }{ \p \nu^{(1)}} \mathcal{S}_{D_2 }   \\
  \frac {\p }{ \p \nu^{(2)}} \mathcal{S}_{D_1 }  &  \Kcal_{D _2}^*
 \end{bmatrix}
 \eeq
(with $I$ being the identity operator), and
$$
\Gvf :=\begin{bmatrix} \one{\Gvf} \\ \two{\Gvf} \end{bmatrix},  \quad
\p h:= \begin{bmatrix}
  \frac {\p h }{\p \nu^{(1)}}  \\
 \nm \frac {\p  h}{\p \nu^{(2)}}
\end{bmatrix} .
$$

If there are $N$ simply connected inclusions, say
$D_1,\ldots,D_N$, then the corresponding NP operator $\Kbb^*$ is
defined by
 \beq
 \Kbb^* :=\begin{bmatrix}
  \Kcal_{D _1}^* & \frac {\p }{ \p \nu^{(1)}} \mathcal{S}_{D_2 } & \cdots & \frac {\p }{ \p \nu^{(1)}} \mathcal{S}_{D_N }  \\
  \frac {\p }{ \p \nu^{(2)}} \mathcal{S}_{D_1 } & \Kcal_{D _2}^* & \cdots   &  \frac {\p }{ \p \nu^{(2)}} \mathcal{S}_{D_N }\\ \vdots & \vdots & \ddots & \vdots\\
  \frac {\p }{ \p \nu^{(N)}} \mathcal{S}_{D_1 } &  \frac {\p }{ \p \nu^{(N)}} \mathcal{S}_{D_2 } & \cdots   &  \Kcal_{D _N}^*
 \end{bmatrix}.
 \eeq

We make note of some important properties of the NP operator
$\Kbb^*$ whose proofs can be found in \cite{ACKLM}. Let $\Hcal:=
L^2(\p D_1) \times L^2(\p D_2)$ and $\Hcal_0:= L^2_0(\p D_1)
\times L^2_0(\p D_2)$. We have
\begin{itemize}
\item $\Kbb^*$ maps $\Hcal$ into itself, and $\Hcal_0$ into itself.
\item For any $\Gl$ with $\Gl \le -1/2$ or $\Gl > 1/2$, $\Gl \Ibb - \Kbb^*$ is invertible on $\Hcal$.
\item $\frac{1}{2}\Ibb -\Kbb^*$ is invertible on $\Hcal_0$.
\item All the eigenvalues of $\Kbb^*$ belong to $(-1/2, 1/2]$.
\end{itemize}

One fact of crucial importance is that $\Kbb^*$ can be
symmetrized. To see this we introduce the operator $\Sbb$ acting
on $\Hcal$ as
 \beq\label{eqSbb}
 \Sbb = \begin{bmatrix}
 \Scal_{D_1} & \Scal_{D_2} \\
 \Scal_{D_1} & \Scal_{D_2}
 \end{bmatrix}.
 \eeq
It is worth making a remark on the operators off diagonal. For
example, $\Scal_{D_2}$ on the top right corner is an operator from
$L^2(\p D_2)$ into $L^2(\p D_1)$. It is proved in \cite{ACKLM}
based on a result in \cite{KPS} that $-\Sbb$ is positive
semi-definite and self-adjoint, $\Sbb\Kbb^*$ is self-adjoint, and
hence there is a self-adjoint operator $\Abb$ on $\Hcal$ such that
 \beq\label{symmetry}
 \sqrt{-\Sbb} \Kbb^* = \Abb \sqrt{-\Sbb}.
 \eeq
In other words, $\Kbb^*$ is self-adjoint with respect to the inner product
 \beq\label{symminner}
 \la \Gvf, \Gy \ra_\Abb := - \la \Sbb[\Gvf], \Gy \ra, \quad \Gvf, \Gy \in \Hcal.
 \eeq

%%%%%%%%%%%%%%%%%%%%%%%%%%%%%%%%%%%%%%%%%%%%
\section{Eigenfunctions of $\Kbb^*$}
%%%%%%%%%%%%%%%%%%%%%%%%%%%%%%%%%%%%%%%%%%%%

We now construct eigenfunctions of $\Kbb^*$ corresponding to
$1/2$. Our construction plays an essential role in understanding
the blow-up of the gradient. We first prove the following lemma.
\begin{lem}\label{eigenvector}
For $i=1,2$, there is a unique solution $v_i$ to
 \begin{equation}\label{def_of_vi}
\quad \left\{
\begin{array}{ll}
\ds \Delta v_i  = 0 \quad\mbox{in } \Rbb^d \setminus \overline{D_1 \cup D_2},\\
\ds v_i = \Gl_j \ (\mbox{constant}) \quad\mbox{on } \p D_j, \ j=1,2,  \\
\int_{\p D_i} \frac {\p v_i}{\p \nu^{(i)}} \Big|_+
d\Gs \neq 0, \ \int_{\p D_j} \frac {\p v_i}{\p \nu^{(j)}} \Big|_+ d\Gs  =0 \quad \mbox{if } j \neq i, \\
\ds v_i(\Bx) - \Scal_{D_i} [1](\Bx)   =O(|\Bx|^{1-d}) \quad \mbox{as } |\Bx|\rightarrow \infty.\end{array}
\right.
\end{equation}
\end{lem}
\pf
We first observe that
 $$
 \int_{\p D_j}  \frac {\p }{\p \nu^{(j)}} \Scal_{D_i} [1]\Big|_- d\Gs  =0, \quad j=1,2.
 $$
Since $\frac 1 2 \Ibb - \Kbb^*$ is invertible on $\Hcal_0$, there exists a unique solution
$(\one{\psi_{i}}, \two{\psi_{i}} )^T \in \Hcal_0$ such that
 \beq\label{psieqn}
 \left( \frac 1 2 \Ibb - \Kbb^* \right) \begin{bmatrix} \one{\psi_{i}} \\ \two{\psi_{i}} \end{bmatrix}
 =  \begin{bmatrix}
\frac {\p }{\p \nu^{(1)}} \Scal_{D_i} [1] \big|_- \\
 \nm \frac {\p }{\p \nu^{(2)}} \Scal_{D_i} [1] \big|_-\end{bmatrix} \in \Hcal_0.
 \eeq
Let $\Gd_{ij}$ be the Kronecker symbol and define
 \beq\label{defGvf}
 \begin{bmatrix} \one{\Gvf_{i}} \\ \two{\Gvf_{i}} \end{bmatrix} := \begin{bmatrix} \one{\psi_{i}} \\ \two{\psi_{i}} \end{bmatrix} +
 \begin{bmatrix} \Gd_{i1} \\ \Gd_{i2} \end{bmatrix}.
 \eeq
Since
 $$
 \begin{bmatrix}
\frac {\p }{\p \nu^{(1)}} \Scal_{D_i} [1] \big|_- \\
 \nm \frac {\p }{\p \nu^{(2)}} \Scal_{D_i} [1] \big|_-\end{bmatrix} = \left( -\frac 1 2 \Ibb + \Kbb^* \right)
 \begin{bmatrix} \Gd_{i1} \\ \Gd_{i2} \end{bmatrix},
 $$
we have
 \beq\label{eigen12}
 \left( \frac 1 2 \Ibb - \Kbb^* \right) \begin{bmatrix} \one{\Gvf_{i}} \\ \two{\Gvf_{i}} \end{bmatrix} =0.
 \eeq
Define
 \beq
 v_i(\Bx): = \Scal_{D_1}[\one{\Gvf_{i}}](\Bx) + \Scal_{D_2}[\two{\Gvf_{i}}](\Bx) ,
 \quad \Bx \in \Rbb^d \setminus (D_1 \cup D_2).
 \eeq

We now show that $v_i$ is the desired function. Because of \eqnref{eigen12}, we have
 $$
 \pd{}{\nu^{(j)}} \big( \Scal_{D_1}[\one{\Gvf_{i}}] + \Scal_{D_2}[\two{\Gvf_{i}}] \big)
 \Big|_{-} =0 \quad\mbox{on } \p D_j,\quad j=1,2,
 $$
and hence, $\Scal_{D_1}[\one{\Gvf_{i}}] +
\Scal_{D_2}[\two{\Gvf_{i}}]$ is constant in $D_1$ and $D_2$. Thus,
$v_i$ is constant on $\p D_1$ and $\p D_2$.

Since $(\one{\psi_{i}},\two{\psi_{i}} )^T \in \Hcal_0$, we have
 $$
 \int_{\p D_j} \pd{}{\nu^{(j)}} \big( \Scal_{D_1}[\one{\psi_{i}}] + \Scal_{D_2}[\two{\psi_{i}}] \big) \Big|_{+} \, d\Gs=0, \quad j=1,2.
 $$
On the other hand, we have
 $$
 \int_{\p D_i} \pd{}{\nu^{(i)}} \Scal_{D_i} [1] \Big|_{+} \,
 d\Gs= \int_{\p D_i} 1+ \pd{}{\nu^{(i)}} \Scal_{D_i} [1] \Big|_{-} \, d\Gs = |\p D_i|,
 $$
and
 $$
 \int_{\p D_j} \pd{}{\nu^{(j)}} \Scal_{D_i} [1] \Big|_{+} \, d\Gs= \int_{\p D_j} \pd{}{\nu^{(j)}}
 \Scal_{D_i} [1] \Big|_{-} \, d\Gs = 0
 $$
if $j \neq i$. Here $|\p D_j|$ denotes the area (or length) of $\p D_j$. Thus we have the third line in \eqnref{def_of_vi}. This completes the proof. \qed

As an immediate consequence we have the following theorem.

\begin{thm}\label{specthm}
The eigenvalue $\frac{1}{2}$ of $\Kbb^*$ has multiplicity $2$.
\end{thm}
\pf
The identity \eqnref{eigen12} shows that $\Gvf_j:= (\one{\Gvf_{j}},
\two{\Gvf_{j}})^T$, $j=1,2$, are two eigenfunctions of $\Kbb^*$ corresponding to $1/2$. We have from \eqnref{defGvf} that
 \beq\label{indep}
 \int_{\p D_i} \Gvf_j^{(i)} d\Gs = |\p D_j| \delta_{ij}.
 \eeq
This implies that $\Gvf_1$ and
$\Gvf_2$ are linearly independent in $\Hcal$. Since $\Hcal_0$ has
codimension $2$ in $\Hcal$ and $\frac{1}{2} \Ibb - \Kbb^*$ is
invertible in $\Hcal_0$ as mentioned before, the multiplicity of
$1/2$ is 2. \qed

Using exactly the same arguments one can generalize Theorem
\ref{specthm} to the case when there are $N$ simply connected
inclusions.
\begin{thm}
If there are $N$ simply connected mutually disjoint inclusions, then the eigenvalue $1/2$ of $\Kbb^*$ has multiplicity $N$.
\end{thm}

%%%%%%%%%%%%%%%%%%%%%%%%%%%%%%%%%%%%%%%%%%%%
\section{Characterization of the gradient blow-up}
%%%%%%%%%%%%%%%%%%%%%%%%%%%%%%%%%%%%%%%%%%%%

Let $\Gvf_j$, $j=1,2$, be the eigenfunctions of $\Kbb^*$ corresponding to $1/2$ introduced in the proof of Theorem \ref{specthm}. Because of \eqnref{indep}, if we define $\Bg$ by
\beq\label{Bg}
\Bg:= \frac{1}{|\p D_1|} \Gvf_1 - \frac{1}{|\p D_2|} \Gvf_2,
\eeq
then we have
\beq\label{plusone2}
\int_{\p D_1} \one{g} d\Gs = 1, \quad
\int_{\p D_2} \two{g} d\Gs =-1 .
\eeq
Define
 \beq\label{defq}
 q(\Bx):= \Scal_{D_1}[\one{g}] (\Bx) + \Scal_{D_2}[\two{g}] (\Bx) , \quad \Bx \in \Rbb^d \setminus (D_1 \cup D_2).
 \eeq
Then, $q$ is the solution to
\begin{equation}\label{qequation}
\quad \left\{
\begin{array}{ll}
\ds \Delta q  = 0 \quad\mbox{in } \Rbb^d \setminus \overline{(D_1 \cup D_2)},\\
\ds q = \mbox{constant} \quad\mbox{on } \p D_i,~i=1,2, \\
\int_{\p D_1} \frac {\p q}{\p \nu^{(1)}} \big|_+ d\Gs =1, \ \ \int_{\p D_2} \frac {\p q}{\p \nu^{(2)}} \big|_+ d\Gs =-1 , \\
\ds q(\Bx) =O(|\Bx|^{1-d}) \quad \mbox{as } |\Bx|\rightarrow \infty, \end{array}
\right.
\end{equation}
and
\beq\label{qequation2}
\begin{bmatrix} \pd{q}{\nu^{(1)}}|_{+} \\ \nm \pd{q}{\nu^{(2)}}|_{+} \end{bmatrix} = \Bg.
\eeq
In fact, since $\Kbb^*[\Bg]=\frac{1}{2}\Bg$, we have
 $$
 \begin{bmatrix} \pd{}{\nu^{(1)}} \big( \Scal_{D_1}[\one{g}] (\Bx) + \Scal_{D_2}[\two{g}] \big)
 \big|_{-} \\
 \pd{}{\nu^{(2)}} \big( \Scal_{D_1}[\one{g}] (\Bx) + \Scal_{D_2}[\two{g}] \big)
 \big|_{-} \end{bmatrix} = \left( -\frac 1 2 \Ibb + \Kbb^* \right) \begin{bmatrix} \one{g} \\ \two{g} \end{bmatrix} =0.
 $$
Thus we have the second line in \eqnref{qequation} and \eqnref{qequation2}. The third line in \eqnref{qequation} follows from \eqnref{qequation2}. Because of \eqnref{plusone2}, we have
\begin{align*}
q(\Bx) &= \int_{\p D_1} \GG(\Bx-\By) \one{g}(\By)~d\sigma (\By) + \int_{\p D_2} \GG(\Bx-\By) \two{g}(\By)~d\sigma (\By) \\
&= \int_{\p D_1} \big[\GG(\Bx-\By)-\GG(\Bx-\By_0) \big] \one{g}(\By)~d\sigma (\By) \\
& \qquad + \int_{\p D_2} [\GG(\Bx-\By)-\GG(\Bx-\By_0) \big] \two{g}(\By)~d\sigma (\By)
\end{align*}
for any fixed $\By_0$. Since
$$
|\GG(\Bx-\By)-\GG(\Bx-\By_0)| \le C |\Bx|^{1-d} \quad\mbox{as } |\Bx| \to \infty,
$$
we have the last line in \eqnref{qequation}.

It is known that if $D_1$ and $D_2$ are disks, then the singular function $q$ is given by \eqnref{singular} and completely characterizes the blow-up of $\nabla u$ (see \cite{KLY}). We have the following theorem as a generalization.
Here and throughout this paper $u|_{\p D_j}$ denotes the (constant) value of $u$ on $\p D_j$.

\begin{prop}\label{thm1}
Assume that $D_1$ and $D_2$ are simply connected domains in $\Rbb^d$, $d=2,3$, with
 $\mathcal{C}^{1,\Ga}$ boundaries for some $\alpha >0$. The solution $u$ to \eqnref{inftycond} can be written as
 \beq\label{decomp}
 u(\Bx)= c_\Ge q(\Bx) + b(\Bx), \quad \Bx \in \Rbb^d \setminus (D_1 \cup
 D_2),
 \eeq
where the constant $c_\Ge$ is given by
 \beq\label{crepre}
 c_\Ge := \frac{u|_{\p D_1} - u|_{\p D_2}}{q|_{\p D_1} - q|_{\p D_2}} =  \frac{\sum_{j=1}^2 \int_{\p D_j} h \pd{q}{\nu^{(j)}} \big|_{+} \, d\Gs}{q|_{\p D_1} - q|_{\p D_2}}  ,
 \eeq
and
 \beq\label{bbound}
 \| \nabla b \|_{L^\infty(\GO)} \le C
 \eeq
for some $C$ independent of $\Ge$. Here $\GO = \GO_1 \setminus (D_1 \cup D_2)$ and $\GO_1$ is an open set containing $\overline{D_1 \cup D_2}$.
\end{prop}

\begin{rem}
The constant $c_\Ge$ in \eqnref{crepre} may depend on $\Ge:=
\mbox{dist}(D_1, D_2)$, but is bounded independently of $\Ge$, and hence the singular function $q$ determines the blow-up of $\nabla u$, if
$D_1$ and $D_2$ are strictly convex and have $\mathcal{C}^{2,\Ga}$
smooth boundaries. This fact will be proved in the next
section.
\end{rem}

\noindent{\sl Proof of Proposition \ref{thm1}}.
Let $u$ be the solution to \eqnref{inftycond} and let
$$
\Bu:=\begin{bmatrix} u|_{\p D_1} \\ u|_{\p D_2} \end{bmatrix}, \quad \Bh:=\begin{bmatrix} h|_{\p D_1} \\ h|_{\p D_2} \end{bmatrix}.
$$
Since $u$ is constant on $\p D_1$ and $\p D_2$, it follows from the third line of \eqnref{qequation} that
$$
u|_{\p D_1} - u|_{\p D_2} = \int_{\p D_1} u \pd{q}{\nu^{(1)}} \big|_{+} \, d\Gs + \int_{\p D_2} u \pd{q}{\nu^{(2)}} \big|_{+} \, d\Gs = \la \Bu, \Bg \ra.
$$
The representation \eqnref{urepres} implies that
$$
\Bu=\Bh + \Sbb[\Gvf]
$$
and $\Gvf \in \Hcal_0$. Therefore we have
$$
\la \Bu, \Bg \ra = \la \Bh + \Sbb[\Gvf], \Bg \ra = \la \Bh, \Bg \ra + \la \Gvf, \Sbb[\Bg] \ra = \la \Bh, \Bg \ra.
$$
So we have the second identity in \eqnref{crepre}.

Let
$$
b(\Bx)= u(\Bx)- c_\Ge q(\Bx)= u(\Bx)- \frac{u|_{\p D_1} - u|_{\p D_2}}{q|_{\p D_1} - q|_{\p D_2}} q(\Bx).
$$
Then we have
$$
b|_{\p D_1} = b|_{\p D_2}.
$$
So, one can show following the same lines of the proof of Theorem 2.1 in \cite{KLY} that \eqnref{bbound} holds.
This completes the proof. \qed

It will be shown that the gradient of the singular function $q$ defined as a single layer potential of an eigenfunction $\Bg$ of $\Kbb^*$ blows up as $\Ge=\mbox{dist}(D_1, D_2) \to 0$.  We now show that another eigenfunction does not contribute to the blow-up. Let $\Bf = (\one{f}, \two{f})^T \in \Hcal$ be an eigenfunction of $\Kbb^*$ orthogonal to $\Bg$ with respect to the inner product \eqnref{symminner}, namely,
 $$
 \la \Sbb[\Bf], \Bg \ra =0.
 $$
Then \eqnref{plusone2} shows that $\one{\Gl}=\two{\Gl}$ where $\Gl^{(j)}$ is the $j$-th component of the (constant) vector $\Sbb[\Bf]$. It implies that the function $v$ defined by
 $$
 v(\Bx):= \Scal_{D_1}[\one{f}] (\Bx) + \Scal_{D_2}[\two{f}] (\Bx)
 $$
is constant on $\p D_j$, $j=1,2$, and satisfies
$$
v|_{\p D_1} = v|_{\p D_2}.
$$
So, $|\nabla v|$ stays bounded regardless of $\Ge$.

%%%%%%%%%%%%%%%%%%%%%%%%%%%%%%%%%%%%%%%%%%%%%%%%%%%
\section{Estimates of the potential difference}
%%%%%%%%%%%%%%%%%%%%%%%%%%%%%%%%%%%%%%%%%%%%%%%%%%%

We assume for the rest of this paper that $D_1$ and $D_2$ are
strictly convex domains in $\Rbb^2$ with $\mathcal{C}^{2,\Ga}$
boundaries for some $\alpha >0$. Let $\Bz_1$ and $\Bz_2$ be points
on $\p D_1$ and $\p D_2$, respectively, such that
 \beq
 |\Bz_1 - \Bz_2| = \mbox{dist} (D_1, D_2)=\Ge.
 \eeq

We prove the following proposition.

\begin{prop}\label{thm2}
Let $u$ be the solution to \eqnref{inftycond}, $c_\Ge$ be the constant defined by \eqnref{crepre}, and $\Gk_j$ be the
curvature of $\p D_j$ at $\Bz_j$ for $j=1,2$. Then $c_\Ge$ is bounded independently of $\Ge$ and
 \beq\label{curform}
 u|_{\p D_1} - u|_{\p D_2} = -\frac{ c_\Ge}{\sqrt 2 \pi} \sqrt{\Gk_1 + \Gk_2} \sqrt\Ge (1+O_\Ga),
 \eeq
where
 \beq\label{curform2}
 O_\Ga =  \left\{
\begin{array}{ll}
\ds  O ({\Ge}^{\alpha/ 2 } )~&\mbox{if }\alpha \in (0,1),\\ \ds  O (|\sqrt{\Ge} \ln\Ge|)  ~&\mbox{if }\alpha =1.\end{array}
\right.
 \eeq
\end{prop}

We prove Proposition \ref{thm2} after a sequence of lemmas.

Let $B_j$ be the osculating disk to $D_j$ at $\Bz_j$ so that its radius $r_j$ equals to $1/\Gk_j$. Let $q_B$ be the singular function associated with $B_1$ and $B_2$, {\it i.e.}, the solution to \eqnref{qequation} with $D_1$ and $D_2$ replaced with $B_1$ and $B_2$. Recall that $q_B$ is given explicitly by
 \beq\label{qfordisk}
 q_B (\Bx) = \frac 1 {2\pi}\left(\ln  |\Bx - \Bp_1| - \ln |\Bx - \Bp_2|
 \right),
 \eeq
where $\Bp_1 \in B_1$ and $\Bp_2 \in B_2$ are the unique fixed points of repeated reflections $R_1 R_2$ and $R_2 R_1$, respectively. We emphasize that $q_B$ is harmonic in $\Rbb^2 \setminus \{\Bp_1, \Bp_2\}$.

For the rest of this paper we assume that $\Bz_1= (-\Ge/2,0)$ and $\Bz_2=(\Ge/2,0)$ after translation and rotation if necessary, so that the centers of $B_1$ and $B_2$ are on the real axis. In this case, $\Bp_1$ and $\Bp_2$ are of the forms $\Bp_j=(p_j, 0)$, $j=1,2$, and it is proved in \cite{Y} that
 \beq\label{pjest}
 p_j=  (-1)^j \sqrt 2 \sqrt {\frac {r_1 r_2}{r_1 + r_2 }} \sqrt {\Ge} + O(\Ge), \quad j=1,2.
 \eeq
It is also proved using \eqnref{pjest} that
 \beq\label{diskest}
 q_B |_{\p B_1} - q_B |_{\p B_2} = -\left( \frac {1} {\sqrt 2 \pi} \sqrt { \frac {r_1 + r_2}{r_1 r_2}}\right) \sqrt{\Ge} + O(\epsilon)~\mbox{as}~\Ge \rightarrow 0.
 \eeq
Using \eqnref{pjest} one can see that \eqnref{qBmax2} holds.

\begin{lem}
There is a constant $C$ independent of $\Ge$ such that
\beq\label{comp1} \left| \pd{q}{\nu^{(j)}} \Big |_{+} (\Bx)\right|
\le C |\nabla q_B(\Bx)| \quad \mbox{for all } \Bx \in \p D_j, \ \
j=1,2. \eeq
\end{lem}
\pf We only prove \eqnref{comp1} for $j=2$ since the case for
$j=1$ can be treated in the exactly same way. We first assume that
$D_2$ is a disk so that $B_2=D_2$. Let $B_1'$ be a disk of radius
$r_1'$ and let $B_1''$ be a disk of radius $r_1''$ ($r_1'$ and
$r_1''$ are independent of $\Ge$) such that $B_1' \subset D_1
\subset B_1''$ and $\Bz_1 \in \p B_1' \cap \p B_1''$. Let $q'$ and
$q''$ be the solution to \eqnref{qequation} with $D_1$ replaced
with $B_1'$ and $B_1''$, respectively. Then the monotonic property
\cite[Lemma 2.4]{LY2} yields
 \beq
 0 \leq -\frac {\p q}{\p\nu^{(2)}} \Big|_+(\Bx) \leq
  - \left(\frac{q''|_{\p B_1''}- q''|_{\p B_2}} {q'|_{\p B_1'}- q'|_{\p B_2}}
  \right) \frac {\p q' }{\p\nu^{(2)}} \Big|_+ (\Bx), \quad \Bx \in \p D_2 .
 \eeq
Because of \eqnref{diskest}, there is a constant $C$ independent of $\Ge$ such that
 \beq
 \frac{q''|_{\p B_1''}- q''|_{\p B_2}} {q'|_{\p B_1'}- q'|_{\p B_2}} \le C,
 \eeq
so we have
 \beq
 \left| \frac {\p q}{\p\nu^{(2)}} \Big|_+(\Bx) \right| \leq C |\nabla q'(\Bx)| , \quad \Bx \in \p D_2.
 \eeq
Note that
 \beq
 q'(\Bx)= \frac 1 {2\pi}\left(\ln  |\Bx - \Bp_1'| - \ln |\Bx - \Bp_2'|
 \right),
 \eeq
where $\Bp_1'$ and $\Bp_2'$ are the fixed points of the repeated reflections with respect to $\p B_1'$ and $\p B_2$.
Using \eqnref{pjest} one can show that
 \beq
 |\nabla q'(\Bx)| \le C |\nabla q_B(\Bx)|, \quad \Bx \in \p D_2
 \eeq
for some constant $C$ independent of $\Ge$. So we have \eqnref{comp1} provided that $D_2$ is a disk.

If $D_2$ is not a disk, we may use a conformal mapping to make it
a disk. In fact, if $\tilde\Psi$ is a conformal mapping from
$\Rbb^2  \cup \{\infty \} \setminus \overline U$ onto $\Rbb^2 \cup \{\infty \} \setminus \overline
{\widetilde{D}_2}$ where
$U$ is the unit disk, then $\tilde\Psi$ can be extended up to $\p
U$ as a $\mathcal{C}^1$ function. Therefore, there are constants
$C_1$ and $C_2$ such that
 \beq\label{conbound}
 C_1 \le |\nabla \tilde\Psi(\Bx)| \le C_2 \quad\mbox{for all } \Bx \in \p U.
 \eeq
Let
 \beq
 \Psi=\tilde\Psi + \Ba_2,
 \eeq
where $\Ba_2$ is defined in \eqnref{ajtrans}. Then $\Psi$ is a
conformal mapping from $\Rbb^2 \setminus \overline U$ onto $\Rbb^2
\setminus \overline {{D}_2}$ and satisfies \eqnref{conbound}.
Moreover, by \cite[Appendix]{Y} and a combination with a  linear fractional transformation, we can also assume that there are two
disks $B_1'$ and $B_1''$ of radii independent of $\Ge$ such that
$B_1' \subset \Psi^{-1}(D_1) \subset B_1''$ and $\p B_1' \cap \p
B_1''$ contains the point on $\p \Psi^{-1}(D_1)$ which is the
closest to $U$. Thus we can apply the same argument as above to
$q\circ \Psi - c_{q}$ to obtain \eqnref{comp1}, where $c_q = \lim_{\Bx \rightarrow \infty} q\circ \Psi (\Bx)$. This completes the proof.
\qed

\begin{lem}\label{lem:potential_diff}
There exists a positive $\Gd_0$
(independent of $\Ge$) such that if $\Bx \in \p D_j$ and $|\Bx
-\Bz_j| \leq \Gd_0$, then
 \beq\label{est1}
 \big| q_B (\Bx) - q_B |_{\p B_j} \big| \leq C \sqrt {\Ge} |\Bx -\Bz_j|^{\alpha},
 \eeq
and
 \beq\label{est2}
 \left| \frac {\p q}{\p\nu^{(j)}}\big|_+ (\Bx)\right| \leq C \frac{\sqrt \Ge}{|\Bx -\Bz_j |^2 + \Ge} .
 \eeq
For any point $\Bx \in \p D_j$ with  $|\Bx-\Bz_j| > \delta_0$,
 \beq\label{est3}
 \big| q_B (\Bx) - q_B |_{\p B_j} \big| \leq C \sqrt {\Ge}
 \eeq
and
 \beq\label{est4}
 \left| \frac {\p q }{ \p \nu^{(j)}}\big|_+ (\Bx) \right| \leq C \sqrt \Ge .
 \eeq
Here, $j=1,2,$ and the constants $C$ are independent of $\Ge$ and
$\Gd_0$.
\end{lem}

\pf
Assume that $j=2$ without loss of generality. There exists $\Gd_0 >0$ (independent of $\Ge$) and functions $x_2, \ x_B: [-\Gd_0, \Gd_0] \rightarrow \mathbb{R}$ such that $x_2(0)=\Ge/2$, $x_2'(0)=0$, $x_B(0)=\Ge/2$, $x_B'(0)=0$, and $\p D_2$ and $\p B_2$ are graphs of $x_2$ and $x_B$ for $|y| \le \Gd_0$, {\it i.e.}, $(x_2(y),y) \in \p D_2$ and $ (x_B(y),y) \in \p B_2$.  We then have
 \beq\label{diffest}
 |x_2(y)| \leq C |y|^2 \quad\mbox{and}\quad |x_2(y) - x_B(y)| \leq C |y|^{2+\alpha}
 \eeq
for some constant $C$ since $D_2$ and $B_2$ are osculating at $\Bz_2$. Since the fixed points $\Bp_1$ and $\Bp_2$ of the repeated reflections are on the $x$-axis, we may write $\Bp_j=(p_j, 0)$, $j=1,2$, and \eqnref{pjest} holds.

If $|y| \leq \Gd_0$, then
\begin{align*}
&\left| q_B (x_2(y), y) - q_B \big|_{\p B_2} \right|\\
&=\frac 1 {2\pi} \big(\ln| (x_2(y)-p_1 ,y)|   - \ln| (x_2(y)-p_2 ,y)| \big)\\
& \quad -\frac 1 {2\pi} \big(\ln| (x_B(y)-p_1 ,y)|   - \ln| (x_B(y)-p_2 ,y)| \big).
\end{align*}
If $0< |y|< \Ge^{\frac 1 {2(2+\alpha)}} \le \Gd_0$, there exists $x_*$ between $x_2(y)$ and $x_B(y)$ such that
\begin{align*}
&\left| q_B (x_2(y), y) - q_B \big|_{\p B_2} \right|\\
&\leq C |x_2(y)- x_B (y)|\left|   \frac {(x_* - p_1)}{ (x_* -p_1)^2 + y^2} - \frac {(x_* - p_2)}{ (x_* -p_2)^2 + y^2}\right|\\
&\leq C |x_2(y)- x_B (y)| |p_1-p_2 |    \frac 1 {y^2}\\
&\leq C  |y|^{2+\alpha} (|y|^{2+\alpha} + \sqrt \Ge)
\frac{ 1} { y^2} \leq C \sqrt {\Ge} |y|^{\alpha},
\end{align*}
where the second to last inequality follows from \eqnref{pjest}.

If $ \Ge^{\frac 1 {2(2+\alpha)}}  \leq |y|\leq \delta_0$, there exists $p_*$ between $p_1$ and $p_2$ such that
\begin{align*}
&\left| q_B (x_2(y), y) - q_B \big|_{\p B_2} \right|\\
&\leq C |p_1- p_2| \left| \frac {(x_2(y) - p_*)}{ (x_2(y) -p_*)^2 + y^2} -  \frac {(x_B(y) - p_*)}{ (x_B (y) -p_*)^2 + y^2}\right|\\
&\leq C |p_1- p_2| |x_2(y) -x_B(y)| \frac 1 {y^2}\\
&\leq C \sqrt {\Ge} (|y|^{2+\alpha} + \sqrt \Ge) \frac{ 1} { y^2} \leq C \sqrt {\Ge} |y|^{\alpha},
\end{align*}
where the the second to last inequality holds because of \eqnref{diffest}.

If $ |(x,y) -(\Ge/2, 0)| > \delta_0$, one can easily see from \eqnref{pjest} that
 $$
 |q_B(x,y)| \le C|p_1-p_2| \le C \sqrt{\Ge},
 $$
and hence we have
$$
\left| q_B (x_2(y), y) - q_B \big|_{\p B_2} \right| \leq \left| q_B (x_2(y), y) \right| + \left| q_B \big|_{\p B_2} \right|\leq  C \sqrt {\Ge}.
$$

Now we estimate $\frac {\p q }{ \p \nu^{(2)}}\big|_+$ on $\p D_2$.
By \eqnref{comp1}, we have
$$
\left| \frac {\p q }{ \p \nu^{(2)}}\big|_+ \right| \leq C |\nabla
q_B|\mbox{ on }\p D_2.
$$
Suppose that $|y| \leq \delta_0$.  Then, we have
\begin{align*}
\left| \frac {\p q }{ \p \nu^{(2)}}\big|_+ (x_2(y),y) \right|
&\leq C \left|\frac {\p  q_B}{\p x}(x_2(y),y)\right|\\
& \leq C \left|   \frac {(x_2(y) - p_1)}{ (x_2(y) -p_1)^2 + y^2} -
\frac {(x_2(y) - p_2)}{ (x_2(y) -p_2)^2 + y^2}\right| .
 \end{align*}
If $|y|< \sqrt {\Ge} $,  then $|x_2(y) -p_j| >  C \sqrt \Ge$ for
$j=1,2$,  and thus we obtain
 $$  |\nabla q (x_2(y),y)| \leq C \frac 1 {\sqrt \Ge}. $$
If $ \sqrt {\Ge} \leq |y| \leq \delta_0$, then it follows that
 $$ |\nabla q (x_2(y),y)| \leq C  |p_1 -p_2| \frac 1 {y^2}\leq C \frac {\sqrt {\Ge}}{y^2}.
$$

For $(x,y)$ with  $|(x,y)-(\Ge/2, 0)| > \delta_0$, we have
$$\left| \nabla q_B (x,y)\right| \leq C \sqrt \Ge ,$$
and \eqnref{est4} follows. This completes the proof. \qed

\begin{lem}\label{prop52}
We have
 \beq\label{520-2}
 q|_{\p D_1} - q|_{\p D_2} = -\frac{1}{\sqrt 2 \pi} \sqrt{\Gk_1 + \Gk_2} \sqrt\Ge  + \left\{
\begin{array}{ll}
\ds  O ({\Ge}^{ {(\alpha +1) }/ 2 } ) \quad &\mbox{if }\alpha \in (0,1),\\
\ds  O (|{\Ge} \ln\Ge|)  \quad &\mbox{if }\alpha =1,\end{array}
\right.
\eeq
as $\Ge \rightarrow 0.$
\end{lem}

\pf
We prove that
 \beq\label{520}
 q|_{\p D_1} - q|_{\p D_2}  = q_B |_{\p B_1} - q_B |_{\p B_2}   + \left\{
\begin{array}{ll}
\ds  O ({\Ge}^{ {(\alpha +1) }/ 2 } ) \quad &\mbox{if }\alpha \in (0,1),\\
\ds  O (|{\Ge} \ln\Ge|)  \quad &\mbox{if }\alpha =1, \end{array}
\right.
\eeq
as $\Ge \rightarrow 0.$ Then \eqnref{520-2} follows from \eqnref{diskest}.

Let
 \beq
 v(\Bx):= q(\Bx)-q_B(\Bx).
 \eeq
Since
 $$
 \int_{\p D_i} \frac {\p q_B}{\p \nu^{(i)}} \Big|_+ d\Gs = \int_{\p B_i} \frac {\p q_B}{\p \nu^{(i)}} \Big|_+ d\Gs, \quad i=1,2,
 $$
the function $v$ satisfies
\beq\label{veqn}
\left\{
\begin{array}{ll}
\ds \Delta v  = 0 \quad\mbox{in } \Rbb^2 \setminus \overline{(D_1 \cup D_2)},\\
\ds v (\Bx) +  q_B (\Bx) - q_B |_{\p B_i}  = \mbox{constant}~\mbox{on } \p D_i, \\
\int_{\p D_i} \frac {\p v}{\p \nu^{(i)}} \Big|_+ d\Gs  = 0 , ~i=1,2,\\
\ds v(\Bx) =O(|\Bx|^{-1}) \quad \mbox{as } |\Bx|\rightarrow
\infty.\end{array} \right. \eeq Then, we have
$$
q|_{\p D_1} - q|_{\p D_2} - \left(q_B |_{\p B_1} - q_B |_{\p B_2} \right)
=   -\sum_{i = 1} ^{2} (-1)^{i}\left( v(\Bx) +  q_B (\Bx) - q_B |_{\p B_i}\right) \Big|_{\p D_i} .
$$
We then obtain from the third line in \eqnref{qequation} and the
second line in \eqnref{veqn} that
$$
q|_{\p D_1} - q|_{\p D_2} - \left(q_B |_{\p B_1} - q_B |_{\p B_2} \right) = \sum_{i = 1} ^{2} \int_{\p D_i }  \left(v +  q_B - q_B
|_{\p B_i} \right)  \frac {\p q}{\p \nu^{(i)}} \Big|_+  d\Gs.
$$
An integration by parts and the third line in \eqnref{veqn} yield
$$
q|_{\p D_1} - q|_{\p D_2} - \left(q_B |_{\p B_1} - q_B |_{\p B_2} \right) = \sum_{i = 1} ^{2} \int_{\p D_i }  \left( q_B - q_B
|_{\p B_i} \right)  \frac {\p q}{\p \nu^{(i)}} \Big|_+  d\Gs.
$$

Let
\begin{align*}
& \left| \int_{\p D_i }  \left( q_B - q_B |_{\p B_i} \right)  \frac {\p q}{\p \nu^{(i)}} \Big|_+  d\Gs \right| \\
&= \left| \int_{|\Bx-\Bz_i|\le \Gd_0} + \int_{|\Bx-\Bz_i| > \Gd_0}
\left( q_B - q_B |_{\p B_i} \right)  \frac {\p q}{\p \nu^{(i)}}
\Big|_+ d\Gs \right| := I_1 + I_2.
\end{align*}
Using \eqnref{est3} and \eqnref{est4} we have
 $$
 |I_2| \le C \Ge.
 $$
To estimate $I_1$, let $N$ be the smallest integer such that $\Gd_0 \le 2^N \sqrt{\Ge}$. We then have from \eqnref{est1} and \eqnref{est2} that
 \begin{align*}
 |I_1| &\le \left| \int_{|\Bx-\Bz_i|\le \sqrt{\Ge}} + \sum_{j=1}^N \int_{2^{j-1} \sqrt{\Ge} < |\Bx-\Bz_i|\le 2^j \sqrt{\Ge}} \left( q_B - q_B |_{\p B_i} \right)
  \frac {\p q}{\p \nu^{(i)}} \Big|_+  d\Gs \right| \\
 & \le C \left[ \int_{|\Bx-\Bz_i|\le \sqrt{\Ge}} |\Bx-\Bz_i| d\Gs + \Ge
 \sum_{j=1}^N \int_{2^{j-1} \sqrt{\Ge} < |\Bx-\Bz_i|\le 2^j \sqrt{\Ge}} \frac{1}{|\Bx-\Bz_i|^{2-\Ga}} d\Gs \right] \\
 &\leq \left\{
\begin{array}{ll}
\ds C_{\Ga}(\Ge + \Ge^{ \frac {1+ \Ga} 2} )\leq 2 C_{\Ga}\Ge^{ \frac{1+\Ga} 2}   \quad &\mbox{if }\Ga \in (0,1),\\
\nm
\ds C(\Ge + \Ge N) \le C \Ge\ln \frac{1}{\Ge} \quad &\mbox{if }\Ga =1.\end{array}
\right. \end{align*}
This completes the proof.
\qed

\medskip
\noindent{\sl Proof of Proposition \ref{thm2}}.
We first prove that $c_\Ge$ is bounded independently of $\Ge$. For that we prove that
\beq\label{hgbound}
|\la \Bh, \Bg \ra| = \left| \sum_{j=1}^2 \int_{\p D_j} h \pd{q}{\nu^{(j)}} \big|_{+} \, d\Gs \right| \le C \sqrt{\Ge}.
\eeq

We still assume that
$\Bz_1= (-\Ge/2,0)$ and $\Bz_2=(\Ge/2,0)$. Pick a point $(c,0) \in D_2$ where $c$ is independent of $\Ge$, and let
\beq
\psi(x,y):= \frac{c^2 y}{(x-c)^2+y^2}.
\eeq
Then, $\psi$ is harmonic except at $(c,0)$ and $\psi(x,y) = O(|(x,y)|^{-1})$ as $|(x,y)| \to \infty$. Since $q$ is constant on $\p D_j$, $j=1,2$, we have by the divergence theorem that
\begin{align}\label{intzero}
\sum_{j=1}^2  \int_{\p D_j} \psi \pd{q}{\nu^{(j)}} \big|_{+}  \, d\Gs= \sum_{j=1}^2  \int_{\p D_j} \pd{\psi}{\nu^{(j)}} \big|_{+} q \, d\Gs = \sum_{j=1}^2  q|_{\p D_j} \int_{\p D_j} \pd{\psi}{\nu^{(j)}} \big|_{+}  \, d\Gs = 0.
\end{align}
Moreover, one can easily see that there is a constant $C>0$ such that
\beq
|\psi(x,y)- y| \le C(x^2 + y^2)\quad \mbox{for all } (x,y)\in \Rbb^2 \setminus (D_1 \cup D_2).
\eeq
Therefore, we have from Taylor's theorem that
\beq\label{taylor}
\left| h (x,y)- h(0,0) - \frac{\p h}{\p y}(0,0) \psi(x,y) \right| \le C (|x| + y^2)
\eeq
for all $(x,y)\in \p D_1 \cup \p D_2$.

Because of the third line in \eqnref{qequation} and \eqnref{intzero}, we have
\begin{align*}
\ds \sum_{j=1}^2 \int_{\p D_j} h \pd{q}{\nu^{(j)}} \big|_{+} \, d\Gs = \sum_{j=1}^2 \int_{\p D_j}\left( h - h(0,0) - \frac{\p h}{\p y}(0,0) \psi \right) \pd{q}{\nu^{(j)}} \big|_{+} \, d\Gs.
\end{align*}
Let
\begin{align*}
\int_{\p D_1} \left( h - h(0,0) - \frac{\p h}{\p y}(0,0) \psi \right) \pd{q}{\nu^{(1)}} \big|_{+} \, d\Gs =
\int_{|\Bx-\Bz_1| \le \Gd_0} + \int_{|\Bx-\Bz_1| > \Gd_0} := I_1+I_2.
\end{align*}
It follows from \eqnref{est2} and \eqnref{taylor} that
$$
|I_1| \le C \int_{|y| \le \Gd_0} \frac{\sqrt \Ge (|x_1(y)| + y^2)}{y^2 + \Ge} \, d\Gs \le C' \int_{|y| \le \Gd_0} \frac{\sqrt \Ge y^2}{y^2 + \Ge} \, d\Gs\le C'' \sqrt{\Ge}.
$$
By \eqnref{est4}, we have
$$
|I_2| \le C \sqrt{\Ge}.
$$
Therefore, we have
$$
\left| \int_{\p D_1} \left( h - h(0,0) - \frac{\p h}{\p y}(0,0) \psi \right) \pd{q}{\nu^{(1)}} \big|_{+} \, d\Gs \right| \le C \sqrt{\Ge}.
$$
Similarly, we can show that
$$
\left| \int_{\p D_2} \left( h - h(0,0) - \frac{\p h}{\p y}(0,0) \psi \right) \pd{q}{\nu^{(2)}} \big|_{+} \, d\Gs \right| \le C \sqrt{\Ge}.
$$
Hence we obtain \eqnref{hgbound}. We now infer from \eqnref{diskest} and Lemma \ref{prop52} that $c_\Ge$ is
bounded regardless of $\Ge$.

Since
 \begin{align*}
 u|_{\p D_1} - u|_{\p D_2} & = c_\Ge (q|_{\p D_1} - q|_{\p D_2})
 \end{align*}
by \eqnref{crepre}, \eqnref{curform} follows from Lemma \ref{prop52}. \qed

%%%%%%%%%%%%%%%%%%%%%%%%%%%%%%%%%%%%%%%%%%%%%%%%%%%
\section{Estimates of the gradient- Proof of Theorem \ref{thm3}}\label{sec6}
%%%%%%%%%%%%%%%%%%%%%%%%%%%%%%%%%%%%%%%%%%%%%%%%%%%

Proposition \ref{thm2} and \eqnref{hgbound} show that
\beq
c_\Ge =- \frac{\sqrt 2 \pi \la \Bh, \Bg \ra}{\sqrt{\Ge(\Gk_1+\Gk_2)}} (1+O_\Ga),
\eeq
where
\beq
O_\Ga =
\left\{
\begin{array}{ll}
\ds
 O ({\Ge}^{ {\Ga}/ 2 } ) ~&\mbox{if }\Ga \in (0,1),\\ \ds O (|\sqrt{\Ge} \ln\Ge|) ~&\mbox{if }\Ga =1. \end{array}
\right.
\eeq
So, Theorem \ref{thm3} is an immediate consequence of Proposition \ref{thm1} and the following proposition.

\begin{prop}\label{thm:gradient_potential}
We have
\beq
q (\Bx) = a_\Ge q_B (\Bx) + v (\Bx), \quad \Bx \in \Rbb^2 \setminus (D_1 \cup D_2),
\eeq
where
\begin{align}\label{65}
a_\Ge :&= {\frac {q|_{\p D_1}- q|_{\p D_2}}{q_B |_{\p B_1}- q_B|_{\p B_2}}}
=1+ \left\{
\begin{array}{ll}
\ds
 O ({\Ge}^{ {\Ga}/ 2 } ) ~&\mbox{if }\Ga \in (0,1),\\ \ds O (|\sqrt{\Ge} \ln\Ge|) ~&\mbox{if }\Ga =1,\end{array}
\right.
\end{align}
and
\beq\label{66}
\norm{\nabla v}_{L^{\infty} (\Rbb^2 \setminus
\overline{(D_1 \cup D_2)})} \leq C
\eeq
for some constant $C$ independent of $\Ge$.
\end{prop}

We first fix notation. We suppose that $\Bz_1= (-\Ge/2,0)$ and $\Bz_2=(\Ge/2,0)$ as before. There exists $\Gd_0 >0$ (independent of $\Ge$) and functions $x_1, \ x_2: [-\delta_0, \delta_0]\rightarrow \mathbb{R}$ such that $x_1(0)=-\Ge/2$, $x_1'(0)=0$, $x_2(0)=\Ge/2$, $x_2'(0)=0$, and $\p D_1$ and $\p D_2$ are graphs of $x_1$ and $x_2$ for $|y| \le \Gd_0$, {\it i.e.}, $(x_1(y),y) \in \p D_1$ and $ (x_2(y),y) \in \p D_2$. Since $D_1$ and $D_2$ are strictly convex, $x_1$ is strictly concave and $x_2$ is strictly convex. For $\Gd \le \Gd_0$, let
$$
\Pi_\Gd: =\{ (x,y) \in \Rbb^2 \setminus (D_1 \cup D_2) ~|~ x_1(y) < x < x_2(y), \  |y| \le \Gd \}.
$$

To prove Proposition \ref{thm:gradient_potential}, we need the
following result whose proof will be given in the last section.

\begin{lem}\label{thm:gradient_general}
If $v$ is a bounded harmonic function in $\mathbb{R}^2 \setminus
\overline{(D_1 \cup D_2)}$ satisfying
 \beq\label{localest1}
 \left|  {\frac {\p^2 v}{\p \tau^2}} (x_i (y),y)\right| \leq M |y|^{\Ga-1}~\mbox{ for } |y| \le \Gd_0,
 \eeq
 \beq\label{tangbound}
 \left \| \frac{\p^2 v}{\p \tau^2} \right\|_{L^{\infty }{((\p D_1 \cup \p D_2) \setminus  \p { \Pi_{\Gd_0}}) } } \leq M
 \eeq
for some constant $M$ (independent of $\Ge$), and
\beq
v(\Bz_1) =\frac{\p}{\p \tau} v (\Bz_1) =v (\Bz_2) =\frac{\p}{\p
\tau} v (\Bz_2) = 0,
\eeq
where $\frac {\p}{\p \tau}$ is the
tangential derivative on $\p D_i$, then there exists a constant
$C$ independent of $\Ge >0$ such that
\beq\label{vgrad}
\norm{\nabla v}_{L^{\infty} (\mathbb{R}^2 \setminus \overline{(D_1 \cup D_2)} )} \leq C.
\eeq
\end{lem}

We also need the following lemma.
\begin{lem}\label{Vero}
There exists a positive constant $M$ independent of $\Ge$ such that
 \beq\label{localest2}
 \left|  {\frac {\p^2 q_B}{\p \tau^2}} (x_i (y),y)\right| \leq M |y|^{\Ga-1}~\mbox{ for } |y| \le \Gd_0 ,
 \eeq and
\beq\label{610}
\left\| \frac{\p^2 q_B}{\p \tau^2} \right\|_{L^{\infty }{(\p D_i \setminus \p \Pi_{\Gd_0 })} } \leq M
\eeq
for $i=1,2$.
\end{lem}

\pf We prove \eqnref{localest2} and \eqnref{610} for $i=2$. We use the same notation as
in the proof of Lemma \ref{lem:potential_diff}:
$\p B_2$ are given by $(x_B(y),y)$ for $|y| \le \Gd_0$. Let $\Bx_2(y) =(x_2(y), y)$ and $\Bx_B (y) =(x_B(y), y)$.

Note that
$$
\frac{\p^2 q_B}{\p \tau^2} = \frac{\p^2}{\p \tau^2} \left(  q_B - q_B (\Bz_2) \right) \approx \frac{d^2}{dy^2} \left( q_B (\Bx_2(y)) - q_B (\Bx_B(y)) \right)
$$
if $|y| \le \Gd_0$. Straightforward computations yield
\begin{align*}
&\frac {d}{dy} \left( q_B (\Bx_2(y)) - q_B (\Bx_B(y)) \right)\\
&=\frac 1 {2\pi} \sum_{i=1}^2 (-1)^{i+1} \left( \frac {(\Bx_2 (y) - \Bp_i)\cdot \Bx_2 '(y)  }{|\Bx_2(y) - \Bp_i|^2} - \frac {(\Bx_B (y) - \Bp_i)\cdot \Bx_B ' (y)  }{|\Bx_B(y) - \Bp_i|^2}  \right),
\end{align*}
and
\begin{align*}
& \frac {d^2 }{dy^2} \left( q_B (\Bx_2(y)) - q_B (\Bx_B(y)) \right)\\
&=\frac 1 {2\pi} \sum_{i=1}^2 (-1)^{i+1} \left( \frac {|\Bx_2 ' (y)|^2  }{|\Bx_2(y) - \Bp_i|^2} - \frac {|\Bx_B ' (y) |^2 }{|\Bx_B(y) - \Bp_i|^2}  \right)\\&~~~+\frac 1 {2\pi} \sum_{i=1}^2 (-1)^{i+1} \left( \frac {(\Bx_2 (y) - \Bp_i)\cdot \Bx_2 '' (y)  }{|\Bx_2(y) - \Bp_i|^2} - \frac {(\Bx_B (y) - \Bp_i)\cdot \Bx_B '' (y)  }{|\Bx_B(y) - \Bp_i|^2}  \right)\\&~~~+\frac 1 {\pi} \sum_{i=1}^2 (-1)^{i} \left( \frac {((\Bx_2 (y) - \Bp_i)\cdot \Bx_2 ' (y)  )^2}{|\Bx_2(y) - \Bp_i|^4} - \frac {((\Bx_B (y) - \Bp_i)\cdot \Bx_B ' (y) )^2 }{|\Bx_B(y) - \Bp_i|^4}  \right) \\
&:= I_1+ I_2 + I_3.
\end{align*}

To estimate $I_1$, $I_2$ and $I_3$, we make some preliminary
computations. Since $B_2$ and $D_2$ are osculating at $\Bz_2$, we
have
\beq\label{611}
|x_2(y) - x_B(y)| \leq C |y|^{2+\Ga}.
\eeq
Since $x_2(0)=x_B(0)=\Ge/2$ and $x_2'(0)=x_B'(0)=0$, we have
\beq
|x_2(y)|+|x_B(y)| \le C (y^2+\Ge), \quad |x_2'(y)|+|x_B'(y)| \le C |y|, \quad |x_2''(y)|+|x_B''(y)| \le C,
\eeq
and
\beq\label{612}
|x_2(y) - p_i| \leq C |y|^2 +
\sqrt{\Ge}, \quad |x_B(y) - p_i| \leq C |y|^2 + \sqrt{\Ge}, \quad
j=1,2.
\eeq
It is worth mentioning that the constant $C$ may
differ at each appearance. We also have
\beq\label{613}
|x_2'(y) - x_B'(y)| \le Cy^{1+\Ga},
\eeq
and
\beq\label{614}
|\Bx_2(y)-\Bp_i|^2 \ge C(y^2+\Ge), \quad |\Bx_B(y)-\Bp_i|^2 \ge C(y^2+\Ge), \quad j=1,2.
\eeq

To estimate $I_1$, we write
\begin{align*}
&\left| \frac {|\Bx_2 ' (y)|^2  }{|\Bx_2(y) - \Bp_i|^2} - \frac {|\Bx_B ' (y) |^2 }{|\Bx_B(y) - \Bp_i|^2}  \right| \\ &\leq |\Bx_2 ' (y)|^2 \left| \frac {1  }{|\Bx_2(y) - \Bp_i|^2} - \frac {1 }{|\Bx_B(y) - \Bp_i|^2}  \right|+ \frac {\left||\Bx_2 ' (y) |^2-|\Bx_B ' (y) |^2 \right|}{|\Bx_B(y) - \Bp_i|^2} .
\end{align*}
Using \eqnref{611}, \eqnref{612} and \eqnref{614} we get
\begin{align}
\left| \frac {1  }{|\Bx_2(y) - \Bp_i|^2} - \frac {1 }{|\Bx_B(y) - \Bp_i|^2}  \right|
& \le \frac{|x_2(y) - x_B(y)| |x_2(y) + x_B(y) -2 p_i| }{|\Bx_2(y) - \Bp_i|^2 |\Bx_B(y) - \Bp_i|^2} \nonumber \\
& \leq C \frac {\sqrt \Ge |y|^{2+\Ga} + |y|^{4+\Ga}}{y^4 +  \Ge^2} \le C\frac {1}{ |y|^{1-\Ga}} . \label{615}
\end{align}
We then use \eqnref{613} to arrive at
\begin{align*}
\frac{\left||\Bx_2 ' (y)|^2-|\Bx_B ' (y) |^2 \right| }{|\Bx_B(y) - \Bp_i|^2}
\le \frac{\left| x_2 ' (y)^2- x_B'(y)^2 \right|}{|\Bx_B(y) - \Bp_i|^2}
\leq C \frac {|y|^{1+\Ga}} {y^2 + \Ge} \leq C  \frac {1}{ |y|^{1-\Ga}}.
\end{align*}
Thus we have
$$
|I_1| \le C|y|^{\Ga-1}.
$$

It follows from \eqnref{615} that
\begin{align*}
&\left| \frac {(\Bx_2 (y) - \Bp_i)\cdot \Bx_2 '' (y)  }{|\Bx_2(y) - \Bp_i|^2} - \frac {(\Bx_B (y) - \Bp_i)\cdot \Bx_B '' (y)  }{|\Bx_B(y) - \Bp_i|^2}  \right|\\
&\leq \left| \frac {\Bx_2 (y) \cdot \Bx_2 '' (y)  }{|\Bx_2(y) - \Bp_i|^2} \right|+\left| \frac {\Bx_B (y) \cdot \Bx_B '' (y)  }{|\Bx_B(y) - \Bp_i|^2}  \right|\\&~~+\left| \frac { \Bp_i \cdot (\Bx_2 '' (y) - \Bx_B '' (y) ) }{|\Bx_2(y) - \Bp_i|^2} \right|+ \left|( \Bp_i\cdot \Bx_B '' (y))\left(\frac {1  }{|\Bx_B(y) - \Bp_i|^2}  - \frac {1  }{|\Bx_2(y) - \Bp_i|^2} \right)\right|\\&\leq C\left( \frac {y^2 +\Ge}{y^2 + \Ge}+ \frac {y^2 +  \Ge}{y^2 + \Ge}+  \frac {\sqrt \Ge |y|^{\Ga}}{y^2 + \Ge} + \sqrt\Ge  \frac {\sqrt \Ge |y|^{2+\Ga} + |y|^{4+\Ga} }{y^4 +  \Ge^2}  \right)\leq C|y|^{\Ga-1},
\end{align*}
and hence
$$
|I_2| \le C|y|^{\Ga-1}.
$$

To estimate $I_3$, we first write
\begin{align*}
&\left| \frac {((\Bx_2 (y) - \Bp_i)\cdot \Bx_2 ' (y)  )^2}{|\Bx_2(y) - \Bp_i|^4} - \frac {((\Bx_B (y) - \Bp_i)\cdot \Bx_B ' (y) )^2 }{|\Bx_B(y) - \Bp_i|^4}  \right|\\
&\leq \left| \frac {(\Bx_2 (y) - \Bp_i)\cdot \Bx_2 ' (y)}{|\Bx_2(y) - \Bp_i|^2}
+ \frac {(\Bx_B (y) - \Bp_i)\cdot \Bx_B ' (y)}{|\Bx_B(y) - \Bp_i|^2}  \right| \\
& \quad \times \left| \frac {(\Bx_2 (y) - \Bp_i)\cdot \Bx_2 ' (y)}{|\Bx_2(y) - \Bp_i|^2}
- \frac {(\Bx_B (y) - \Bp_i)\cdot \Bx_B ' (y)}{|\Bx_B(y) - \Bp_i|^2}  \right|.
\end{align*}
One can easily see from \eqnref{614} that
\beq\label{plus}
\left| \frac {(\Bx_2 (y) - \Bp_i)\cdot \Bx_2 ' (y)}{|\Bx_2(y) - \Bp_i|^2}
+ \frac {(\Bx_B (y) - \Bp_i)\cdot \Bx_B ' (y)}{|\Bx_B(y) - \Bp_i|^2}  \right|
\leq \frac {C |y|}{y^2 + \Ge}.
\eeq
Note that
\begin{align*}
& \frac {(\Bx_2 (y) - \Bp_i)\cdot \Bx_2'(y)}{|\Bx_2(y) - \Bp_i|^2}
- \frac {(\Bx_B (y) - \Bp_i)\cdot \Bx_B'(y)}{|\Bx_B(y) - \Bp_i|^2}  \\
& = \frac {(x_2 (y) - p_i) x_2'(y)}{|\Bx_2(y) - \Bp_i|^2} - \frac {(x_B (y) - p_i)x_B ' (y)}{|\Bx_B(y) - \Bp_i|^2}
 + \frac {y}{|\Bx_2(y) - \Bp_i|^2} - \frac {y}{|\Bx_B(y) - \Bp_i|^2}  \\
& = \frac {x_2 (y) x_2'(y)}{|\Bx_2(y) - \Bp_i|^2} - \frac {x_B (y) x_B '(y)}{|\Bx_B(y) - \Bp_i|^2}
  -\frac {p_i(x_2'(y) - x_B'(y))}{|\Bx_2(y) - \Bp_i|^2} \\
& \quad   + (p_i x_B'(y)-y) \left( \frac {1 }{|\Bx_B(y) - \Bp_i|^2} - \frac {1 }{|\Bx_2(y) - \Bp_i|^2} \right).
\end{align*}
We estimate each term using \eqnref{611}-\eqnref{614} to have
\begin{align*}
& \left| \frac {(\Bx_2 (y) - \Bp_i)\cdot \Bx_2'(y)}{|\Bx_2(y) - \Bp_i|^2}
- \frac {(\Bx_B (y) - \Bp_i)\cdot \Bx_B'(y)}{|\Bx_B(y) - \Bp_i|^2} \right|  \\
& \leq C \left[ \frac{|y|(|y|^2+\Ge)}{y^2 + \Ge} + \frac{|y|(|y|^2+\Ge)}{y^2 + \Ge}
+ \frac{\sqrt{\Ge} |y|^{1+\Ga}}{y^2 + \Ge}
+ (\sqrt{\Ge}+1) |y| \frac {\sqrt \Ge |y|^{2+\Ga} +|y|^{2+\Ga}(y^2+\Ge)}{y^4 + \Ge^2} \right] \\
& \leq C |y|^{\Ga} .
\end{align*}
Combining this estimates with \eqnref{plus} we obtain
$$
|I_3| \le C|y|^{\Ga-1}.
$$

If $\Bx \in \p D_2$ satisfies $|\Bx-(\Ge/2, 0)| > \delta_0$, it can be easily seen that
$$
\left| \frac{\p^2 v}{\p\tau^2} (\Bx) \right| \le M.
$$
The proof is complete.
\qed

\medskip
\noindent{\sl Proof of Proposition \ref{thm:gradient_potential}}.
Note that the second identity in \eqnref{65} follows from Lemma \ref{prop52}, and it shows in particular that $a_\Ge$ is bounded regardless of $\Ge$.

Let
\beq\label{vqa}
v(\Bx) = q(\Bx) - a_\Ge q_B(\Bx) ,
\eeq
and
\beq
w(\Bx):= \frac{1}{a_\Ge} (v(\Bx)-v(\Bz_2)).
\eeq
Then one can see from the definition \eqnref{65} of $a_\Ge$ and \eqnref{vqa} that
\beq
w(\Bx)= q_B(\Bz_i) - q_B(\Bx), \quad \Bx \in \p D_i, \ i=1,2.
\eeq
Since $D_i$ and $B_i$ are osculating at $\Bz_i$, we have in particular
$$
w (\Bz_1)=\frac{\p w}{\p \tau} (\Bz_1) =w (\Bz_2)=\frac{\p w}{\p \tau} (\Bz_2)=0.
$$
It follows from Lemma \ref{Vero} and Lemma \ref{thm:gradient_general} that
$$
\norm {\nabla w}_{L^{\infty }{(\mathbb{R}^2 \setminus \overline{(D_1 \cup D_2)} )}} \leq C.
$$
Since $a_\Ge$ is bounded, we obtain \eqnref{66}. This completes the proof.
\qed

%%%%%%%%%%%%%%%%%%%%%%%%%%%%%%%%%%%%%%%%%%%%%%%
\section{The insulating case}
%%%%%%%%%%%%%%%%%%%%%%%%%%%%%%%%%%%%%%%%%%%%%%%

In this section we deal with the case when the inclusions are
insulating, namely, the problem \eqnref{zeocond}. We closely
follow the argument provided in \cite{KLY}.

Let $h^\perp$ be a harmonic conjugate of $h$, {\it i.e.}, $h+i
h^\perp$ is analytic. Let $u^\perp$ be the solution to
\eqnref{inftycond} with $h^\perp$ in place of $h$. Then the
solution $u$ to \eqnref{zeocond} is a harmonic conjugate of
$u^\perp$ in $\Rbb^2\setminus\overline{D_1\cup D_2}$. By Theorem
\ref{thm3}, we have
\beq\label{charcaterize2}
u^\perp(\Bx) = -\frac{\sqrt 2 \pi \la \Bh^\perp, \Bg \ra}{\sqrt{\Ge(\Gk_1+\Gk_2)}}
\Gb_\Ge q_B (\Bx) + r (\Bx), \quad \Bx \in \Rbb^2 \setminus (D_1
\cup D_2),
\eeq
where $\Gb_\Ge$ is a constant of the form
\beq\label{defGb}
\quad\Gb_\Ge =  \left\{
\begin{array}{ll}
\ds
1+ O ({\Ge}^{ {\Ga}/ 2 } ) \quad &\mbox{if }\Ga \in (0,1),\\
\ds 1+ O (|\sqrt{\Ge} \ln\Ge|) \quad &\mbox{if }\Ga =1,
\end{array}
\right.
\quad\mbox{as } \Ge \to 0.
\eeq

Let $\arg:\mathbb{R}^2\setminus \{(0,0)\} \rightarrow [-\pi, \pi)$ be the argument function with a branch cut along the negative real axis, where $\Bx = (x_1,x_2)$ is identified with $x_1 + i x_2$. Define
 \beq\label{hbot}
 q_B^\perp(\Bx) = \frac 1 {2\pi} \Bigr(\arg (\Bx - \Bp_1)-\arg (\Bx - \Bp_2) - \arg (\Bx - \Bc_1)+ \arg (\Bx - \Bc_2)\Bigr),
 \eeq
where $\Bc_j$ is the center of $B_j$, $j=1,2$.  Note that
$q_B^\perp$ is a harmonic function well defined in $\Rbb^2
\setminus \overline{(B_1\cup B_2)}$ since the jump discontinuity
of the argument function across the branch cut is cancelled out
owing to the fact $\Bp_j, \Bc_j\in B_j$, $j=1,2$. Since $\arg (\Bx
- \Bp_1)-\arg (\Bx - \Bp_2)$ is a harmonic conjugate of $q_B$
except on the branch cut and $|\nabla (\arg (\Bx - \Bc_1) - \arg
(\Bx - \Bc_2))|$ is bounded independently of $\Ge$, we arrive at
the following result.
\begin{thm}\label{thm4}
Let $u$ be the solution to \eqnref{zeocond}. Under the same hypothesis as in Theorem \ref{thm3}, we have
\beq
u(\Bx) = -\frac{\sqrt 2 \pi \la \Bh^\perp, \Bg \ra}{\sqrt{\Ge(\Gk_1+\Gk_2)}}
\Gb_\Ge q_B^\perp (\Bx) + r (\Bx), \quad \Bx \in \GO \setminus (D_1 \cup D_2),
\eeq
where $\Gb_\Ge$ is a constant of the form \eqnref{defGb} and
\beq
\norm{\nabla r}_{L^{\infty} (\GO)} \leq C
\eeq
for some $C$ independent of $\Ge$.
\end{thm}

%%%%%%%%%%%%%%%%%%%%%%%%%%%%%%%%%%%%%%%%%%%%%%%
\section{Boundary value problems}
%%%%%%%%%%%%%%%%%%%%%%%%%%%%%%%%%%%%%%%%%%%%%%%

Let $\GO$ be a bounded domain in $\Rbb^2$ with $\mathcal{C}^2$
boundary. Suppose that $\GO$ contains two perfectly conducting
inclusions $D_j$, $j=1, 2$, which have $\mathcal{C}^{2,\Ga}$
boundaries for some $\Ga>0$. We assume that the inclusions are
away from $\p\GO$, namely, there is a constant $c_0$ such that
 \beq\label{away}
 \mbox{dist}(D_j, \p\GO) \ge c_0, \quad j=1, 2.
 \eeq
In this section we consider the following boundary value problem:
 \beq\label{bvp}
 \quad \left\{
 \begin{array}{ll}
 \ds  \Delta u=0 \quad \mbox{in } \GO \setminus \overline{D_1 \cup D_2},\\
 \ds u= \mbox{constant on } \p D_j, \ j=1, 2, \\
 \nm
 \ds \int_{\p D_j} \pd{u}{\nu^{(j)}} \,ds =0 , \ j=1,2,
 \end{array}
 \right.
 \eeq
with either Dirichlet or Neumann boundary conditions on $\p\GO$.

Define a harmonic function $h$ in $\GO$ by \eqnref{hbvp}. This $h$
plays the role of $h$ in the free space problem \eqnref{inftycond}
and we obtain the following result using exactly the same
arguments as those in Theorem \ref{thm3}.

\begin{thm}\label{thm5}
Let $u$ be the solution to \eqnref{bvp} and let $h$ be the
function defined by \eqnref{hbvp}. Then, under the same hypothesis as in Theorem \ref{thm3}, we have
\beq
u(\Bx) = -\frac{\sqrt 2 \pi \la \Bh, \Bg \ra}{\sqrt{\Ge(\Gk_1+\Gk_2)}} \Ga_\Ge q_B
(\Bx) + r (\Bx), \quad \Bx \in \GO \setminus (D_1 \cup D_2),
\eeq
where $\Ga_\Ge$ is a constant of the form
\beq
\Ga_\Ge =  \left\{
\begin{array}{ll}
\ds 1+ O ({\Ge}^{ {\Ga}/ 2 } ) \quad &\mbox{if }\Ga \in (0,1),\\
\ds 1+ O (|\sqrt{\Ge} \ln\Ge|) \quad &\mbox{if }\Ga =1,\end{array}
\right.
\quad\mbox{as } \Ge \to 0,
\eeq
and
\beq
\norm{\nabla r}_{L^{\infty} (\GO)} \leq C
\eeq
for some $C$ independent of $\Ge$.
\end{thm}

One can also obtain a similar result for the insulating boundary
value problem.

%%%%%%%%%%%%%%%%%%%%%%%%%%%%%%%%%%%%%%%%%%%%%%%%%%%%%%%%%
\section{Proof of Lemma \ref{thm:gradient_general}}
\label{sect9}
%%%%%%%%%%%%%%%%%%%%%%%%%%%%%%%%%%%%%%%%%%%%%%%%%%%%%%%%%

We use the same notation as in Section \ref{sec6}.

Let $\chi_2(\Bx)$ be a smooth function on $\p D_2$ such that
\beq
\quad  \left\{
\begin{array}{ll}
\ds \chi_2 (x_2(y),y) = 1      &~\mbox{for} ~|y|\leq {\frac {\delta_0}{2}}, \\
\ds 0 \leq \chi_2 (x_2(y),y) \leq 1  & ~\mbox{for} ~{\frac {\delta_0}{2}}\leq |y|\leq  \delta_0,\\
\ds \chi_2(\Bx) = 0 &~\mbox{otherwise}.\end{array}
\right.
\eeq

As before we denote the tangential derivative on $\p D_j$ by $\frac{\p}{\p\tau}$. Let $g_{2+}(\Bx)$ and $g_{2-}(\Bx)$ be non-negative functions defined for $\Bx=(x_2(y), y) \in \p D_2$, $|y| \le \Gd_0$, such that
\beq
g_{2+}(\Bz_2)= \frac{\p g_{2+}}{\p\tau}(\Bz_2)= g_{2-}(\Bz_2)= \frac{\p g_{2-}}{\p \tau}(\Bz_2)=0
\eeq
and
\beq
\frac{\p^2 g_{2+}}{\p \tau^2}(\Bx) = \max \left\{ \frac {\p^2 v (\Bx)}{\p \tau^2}, 0 \right\}, \quad
\frac{\p^2 g_{2-}}{\p \tau^2}(\Bx) = \max \left\{ -\frac {\p^2 v (\Bx)}{\p \tau^2}, 0 \right\}.
\eeq
Then, $g_{2+}$ and $g_{2-}$ satisfy
\beq
g_{2+}(x_2(y), y) - g_{2-} (x_2(y), y) = v(x_2(y), y), \quad |y| \le \Gd_0.
\eeq

Let $V_{2+}$, $V_{2-}$, $\widetilde V_{2+}$ and $\widetilde
V_{2-}$ be bounded harmonic functions in $\mathbb{R}^2 \setminus
\overline{D_2}$ which satisfy the following Dirichlet boundary
conditions on $\p D_2$: \beq \left\{
\begin{array}{l}
\ds V_{2+}  =  \chi_2 g_{2+} \\
\ds V_{2-}  = \chi_2 g_{2-} \\
\ds \widetilde V_{2+} =  \left(1- \chi_2 \right) \max \{ v,0\} \\
\ds \widetilde V_{2-} =  - \left(1- \chi_2 \right) \min \{ v ,0\}
\end{array}
\right.
\quad\mbox{on } \p D_2.
\eeq
Then by the maximum principle, $V_{2+}$, $V_{2-}$, $\widetilde V_{2+}$ and $\widetilde V_{2-}$
are non-negative and satisfy
\beq
V_{2+} - V_{2-} + \widetilde V_{2+} - \widetilde V_{2-} = v \quad \mbox{on } \p D_2.
\eeq

Let $v_{2+}$, $v_{2-}$, $\widetilde v_{2+}$ and $\widetilde
v_{2-}$ be bounded harmonic functions in $\mathbb{R}^2 \setminus
\overline {(D_1 \cup D_2)}$ which satisfy the following Dirichlet
conditions on $\p D_1$ and $\p D_2$: \beq v_{2+} = v_{2-}
=\widetilde v_{2+} =\widetilde v_{2-}= 0 \quad\mbox{on } \p D_1 ,
\eeq and \beq \quad  \left\{
\begin{array}{l}
\ds v_{2+} =  V_{2+} \\
\ds v_{2-} =  V_{2-} \\
\ds \widetilde v_{2+ } =\widetilde V_{2+ } \\
\ds \widetilde v_{2- } =\widetilde V_{2- }
\end{array}
\right.
\quad\mbox{on } \p D_2.
\eeq
Since $0 \leq v_{2\pm} \leq V_{2\pm}$ and $0 \le \widetilde v_{2\pm} \le \widetilde V_{2\pm}$ on $\p D_1$ and $\p D_2$, we have
\beq \label{V2+}
0 \leq v_{2\pm} \leq V_{2\pm} \mbox{ and } 0 \le \widetilde v_{2\pm} \le \widetilde V_{2\pm} \quad\mbox{in}~ \mathbb{R}^2 \setminus \overline {(D_1 \cup D_2)} .
\eeq

Let \beq w_2:= v_{2+} - v_{2-}+\widetilde v_{2+} - \widetilde
v_{2-}. \eeq Then $w_2$ is a bounded harmonic function in
$\mathbb{R}^2 \setminus \overline {(D_1 \cup D_2)}$ and is such
that \beq w_2 = 0   ~\mbox{on} ~\p D_1, \quad w_2 = v  ~\mbox{on}
~\p D_2. \eeq

In the same way, non-negative bounded harmonic functions $v_{1+}$, $v_{1-}$, $\widetilde v_{1+}$ and $\widetilde v_{1-}$ in $\mathbb{R}^2 \setminus \overline {(D_1 \cup D_2)} $ can be defined so that $w_1:= v_{1+}- v_{1-}+\widetilde v_{1+}- \widetilde v_{1-}$ satisfies
\beq
w_1 = v   ~\mbox{on} ~\p D_1, \quad w_1 = 0  ~\mbox{on} ~\p D_2.
\eeq
Then, we have from the uniqueness of the Dirichlet problem that
\beq\label{sumV}
v= w_1+w_2 \quad \mbox{in} ~\mathbb{R}^2 \setminus \overline {(D_1 \cup D_2)} .
\eeq

We first estimate $|\nabla v_{2+}|$. Thanks to \eqnref{localest1} and \eqnref{tangbound}, we have

\beq\label{boundary}
 \left|  {\frac {\p^2 V_{2+}}{\p \tau^2}} (x_2 (y),y)\right| \leq C |y|^{\Ga-1}~\mbox{ for } |y| \le \Gd_0 ,
 \eeq
 and
 \beq \label{boundary2} \left \| \frac{\p^2
V_{2+}}{\p \tau^2} \right\|_{L^{\infty }{(\p D_2 \setminus  \p { \Pi_{\Gd_0}}) } } \leq C.
\eeq

Since $V_{2+} (\Bz_2) =\frac {\p }{\p \tau }V_{2+}(\Bz_2) = 0$, we have $\| V_{2+} \|_{\mathcal{C}^{1,\Ga}(\p D_2)} \le C$. A standard regularity estimate for harmonic functions yields, in particular,
\beq\label{lem54-2}
\| V_{2+} \|_{\mathcal{C}^{1} (\mathbb{R}^2 \setminus D_2 )} \leq C.
\eeq
If $(x,y) \in \Pi_{\Gd_0}$, then
$$
|(x,y)- (x(y),y)| \leq  C (y^2 + \epsilon).
$$
Thus we obtain from \eqref{lem54-2} and the mean value theorem that
\beq\label{lem54-1}
0 \leq V_{2 +} (x,y) \leq C ( y^2 + \epsilon ), \quad (x,y) \in \Pi_{\Gd_0}.
\eeq
It then follows from \eqnref{V2+} that
\beq\label{v2+bound}
0 \le v_{2+}(x,y) \leq C ( y^2 + \epsilon ), \quad (x,y) \in \Pi_{\Gd_0}.
\eeq

Let $d(\Bx):= \mbox{dist}(\Bx, \p D_2)$ and $B_{r}(\Bx)$ be the disk of radius $r$ centered at $\Bx$. Since $v_{2+}=0$ on $\p D_1$ and $\p D_1$
is $\mathcal{C}^{2,\Ga}$, by a standard elliptic regularity estimate we have
\beq\label{ellreg}
|\nabla v_{2+}(\Bx)| \le \frac{C \| v_{2+} \|_{L^\infty(B_{d(\Bx)}(\Bx) \cap (\mathbb{R}^2 \setminus (D_1 \cup D_2)))}}{d(\Bx)}.
\eeq

If $\Bx=(x,y) \in \Pi_{\Gd_0/2}$ and $x<0$, then there are $c_1$ and $c_2$ such that
 $$
 c_1(y^2+\Ge) \le d(\Bx) \le c_2(y^2+\Ge).
 $$
If $(x', y') \in B_{d(\Bx)}(\Bx) \cap (\mathbb{R}^2 \setminus (D_1 \cup D_2))$, then $|y'| \le |y|+ c_2 (y^2+\Ge)$, thus we have
$$
|v_{2+} (x',y')| \le C(y'^2+\Ge) \le C'(y^2+\Ge).
$$
We then get from \eqnref{ellreg} that
$$
|\nabla v_{2+}(\Bx)| \le C.
$$

If $\Bx=(x,y) \in \Pi_{\Gd_0/2}$ and $x \ge 0$, then we can apply the same argument using the fact that $V_{2+}-v_{2+}=0$ on $\p D_2$ to obtain
$$
|\nabla (V_{2+} - v_{2+})(\Bx)| \le C.
$$
We then obtain using \eqnref{lem54-2} that
$$
|\nabla v_{2+}(\Bx)| \le C.
$$
So, we have
\beq
|\nabla v_{2+}(\Bx)| \le C, \quad \Bx \in \Pi_{\Gd_0/2}.
\eeq

We now estimate $|\nabla \widetilde v_{2+}|$. Since $\| v \|_{\mathcal{C}^{1,\Ga}(\p D_2)} \le C$, we first obtain from the maximum principle that
\beq\label{tildebound}
\| \widetilde V_{2+} \|_{L^\infty(\mathbb{R}^2 \setminus D_2)} \le C
\eeq
for some $C$ independent of $\Ge$. Since
\beq\label{bdryzero}
\widetilde V_{2+}(x_2(y), y)=0 \quad\mbox{if } |y| \le \Gd_0/2,
\eeq
we have
\beq
|\nabla \widetilde V_{2+} (\Bx)| \le C
\eeq
for all $\Bx \in \mathbb{R}^2 \setminus D_2$ satisfying
\beq
\mbox{dist}(\Bx, \p D_2 \setminus \p \Pi_{\Gd_0/2}) \ge \Gd_0/4.
\eeq
In particular, we have
\beq\label{narrowboundtilde}
\sup_{\Bx \in \Pi_{\Gd_0/4}} |\nabla \widetilde V_{2+} (\Bx)| \le C.
\eeq

It follows from \eqnref{bdryzero} and \eqnref{narrowboundtilde} that
$$
0 \leq \widetilde V_{2 +} (x,y) \leq C ( y^2 + \epsilon ), \quad (x,y) \in \Pi_{\Gd_0/4},
$$
and from \eqnref{V2+} that
$$
0 \le \widetilde v_{2+}(x,y) \leq C ( y^2 + \epsilon ), \quad (x,y) \in \Pi_{\Gd_0/4}.
$$
Since $\widetilde v_{2+}(x_1(y),y)= \widetilde v_{2+}(x_2(y),y)=0$ if $|y|\le \Gd_0/2$,  we may apply the same argument as for $v_{2+}$ to obtain
\beq
|\nabla \widetilde v_{2+}(\Bx)| \le C, \quad \Bx \in \Pi_{\Gd_0/8}.
\eeq

In exactly the  same way, one can show that
\beq
|\nabla v_{2-}(\Bx)| + |\nabla \widetilde v_{2-}(\Bx)| \le C, \quad \Bx
\in \Pi_{\Gd_0/8}.
\eeq
Therefore, we have
\beq\label{928}
\sup_{\Bx \in \Pi_{\Gd_0/8}} |\nabla w_2(\Bx)| \le C .
\eeq

If $\Bx \in \p D_1 \setminus \p \Pi_{\Gd_0/8}$, then $d(\Bx) \ge
C$ for some $C$ independent of $\Ge$. Since $w_2=0$ on $\p D_1$ and $w_2$ is bounded, we obtain
\beq\label{929} \sup_{\Bx \in \p D_1 \setminus \p \Pi_{\Gd_0/8}}
|\nabla w_2 (\Bx)| \le C .\eeq Let \beq V= V_{2+} - V_{2-} +
\widetilde V_{2+} - \widetilde V_{2-} \quad\mbox{in } \Rbb^2
\setminus D_2. \eeq Since $V-w_2=0$ on $\p D_2$, it follows that
$$
|\nabla (V-w_2) (\Bx)| \le C
$$
for $\Bx \in \p D_2 \setminus \p \Pi_{\Gd_0/8}$. Since $\| V \|_{C^{1,\Ga}(\p D_2)}$ is bounded, $\| \nabla V \|_{L^\infty(\Rbb^2 \setminus D_2)}$ is bounded, so we have
\beq\label{931}
\sup_{\Bx \in \p D_2 \setminus \p \Pi_{\Gd_0/8}} |\nabla w_2 (\Bx)| \le C.
\eeq

Inequalities \eqnref{928}, \eqnref{929} and \eqnref{931} imply that
$$
\sup_{\Bx \in \p ((\mathbb{R}^2 \setminus (D_1 \cup D_2)) \setminus \Pi_{\Gd_0/8})} |\nabla w_2 (\Bx)| \le C.
$$
We then obtain from the maximum principle that
$$
\sup_{\Bx \in (\mathbb{R}^2 \setminus (D_1 \cup D_2)) \setminus \Pi_{\Gd_0/8}} |\nabla w_2 (\Bx)| \le C.
$$
Combining this with \eqnref{928}, we readily get \beq \sup_{\Bx
\in \mathbb{R}^2 \setminus (D_1 \cup D_2)} |\nabla w_2 (\Bx)| \le
C. \eeq

One can show in exactly  the same way (by switching the roles of
$D_1$ and $D_2$) that \beq \sup_{\Bx \in \mathbb{R}^2 \setminus
(D_1 \cup D_2)} |\nabla w_1 (\Bx)| \le C. \eeq Thus we have
\eqnref{vgrad} and the proof is complete. \qed

%%%%%%%%%%%%%%%%%%%%%%%%%%%%%%%%%%%%%%%%%%%%%%%%%%

\end{document}